\documentclass[aip,
graphicx,
amsmath,
amssymb,
preprint
]{revtex4-1}

\usepackage{float}
\usepackage{subcaption}
\usepackage{mathtools}

\draft % marks overfull lines with a black rule on the right

\begin{document}

% Use the \preprint command to place your local institutional report number 
% on the title page in preprint mode.
% Multiple \preprint commands are allowed.
%\preprint{}

\title{Mixed global dynamics of forced vibro-impact oscillator with Coulomb friction} %Title of paper

% repeat the \author .. \affiliation  etc. as needed
% \email, \thanks, \homepage, \altaffiliation all apply to the current author.
% Explanatory text should go in the []'s, 
% actual e-mail address or url should go in the {}'s for \email and \homepage.
% Please use the appropriate macro for the type of information

% \affiliation command applies to all authors since the last \affiliation command. 
% The \affiliation command should follow the other information.

\author{Oleg Gendelman}
\email[]{ovgend@technion.ac.il}
%\homepage[]{Your web page}
%\thanks{}
%\altaffiliation{}
\affiliation{Faculty of Mechanical Engineering, Technion -- Israel Institute of Technology}

\author{Pavel Kravetc}
\email[]{pavel.kravetc@utdallas.edu}

\author{Dmitrii Rachinskii}
\email[]{Dmitry.Rachinskiy@utdallas.edu}
\affiliation{Department of Mathematical Sciences, The University of Texas at Dallas}

% Collaboration name, if desired (requires use of superscriptaddress option in \documentclass). 
% \noaffiliation is required (may also be used with the \author command).
%\collaboration{}
%\noaffiliation

\date{\today}

\begin{abstract}
The paper revisits a well-known model of forced vibro-impact oscillator with Amonton-Coulomb friction. In vast majority of the existing studies, this model included also viscous friction, and its global dynamics in the state space is governed by periodic, quasiperiodic or chaotic attractors. We demonstrate that removal of the viscous friction leads to qualitative modification of the global dynamics. Namely, the state space is divided into the regions with ``regular" attraction to the aforementioned special solutions, and the regions with profoundly Hamiltonian dynamics. The latter regions contain structures typical for forced Hamiltonian systems: stability islands, extended non-attractive chaotic regions etc.  We prove that such local Hamiltonian behavior should occur for phase trajectories with non-vanishing velocity. Stability analysis for the periodic orbits confirms the above statement. It is demonstrated that similar mixed global dynamics can be observed in broader class of models.
\end{abstract}

% insert suggested PACS numbers in braces on next line

\maketitle %\maketitle must follow title, authors, abstract and \pacs

\begin{quotation}
It is a well-known fact that there is a drastic qualitative difference between Hamiltonian and dissipative dynamics. Hamiltonian systems are governed by Liouville's theorem which imposes strict restrictions on possible dynamics. In particular, it prohibits the existence of asymptotically stable solutions. Dissipative systems, by contrast, usually demonstrate eventual attraction of all phase trajectories to $\omega$-limit sets such as fixed points, periodic orbits, quasiperiodic orbits, or strange attractors. Although it is commonly believed that conservative and dissipative dynamics can not be observed together, it is not true. In this paper, we present a class of dissipative mechanical systems which demonstrate coexistence of both attracting and Hamiltonian-like behaviors.
\end{quotation}

% Body of paper goes here. Use proper sectioning commands. 
% References should be done using the \cite, \ref, and \label commands
\section{\label{sec:intro}Introduction}
Impact oscillators of various types attract a lot of attention for decades due to their importance in various technical applications \cite{babitski, fidlin, kobrinski}. At the same time, dynamics of discontinuous systems involving impacts raises fundamental mathematical questions that inspire major research efforts over decades \cite{bernardo2008piecewise, filippov2013differential, Fredriksson315, ivanov1994}. The impact interaction represents arguably the strongest possible nonlinearity of interaction. On the other hand, unlike other systems with strong nonlinearity, it is ubiquitous, easily realizable and used for modeling multiple physical phenomena \cite{manevichgendelman, manevitch2011tractable}. It is not surprising that the behavior of systems with impacts attracted a lot of attention since the very first s.pdf of nonlinear dynamics \cite{shawholmes, lirandmoon, bapat}.

Amonton-Coulomb dry friction is even more ubiquitous and widely studied dynamical phenomenon, also characterized by a peculiar discontinuity and raising quite a few problems in mathematics and dynamics \cite{feigin, shaw1986}. Dynamics of combined system, {\it i.e.} forced vibro-impact oscillator with dry friction presumably has been first addressed in paper \cite{cone1995}. This work addressed periodic responses of the system and their bifurcations in the space of parameters. This line of research has been continued in a number of studies \cite{virgin1999, Blazejczyk1996, zhang2017, fan2018}, including recent work on systems with multiple degrees of freedom.

It is a common knowledge that Hamiltonian systems, even with time-dependent Hamiltonian, obey Liouville's theorem of the phase volume conservation \cite{landaulifshitz, arnold}. This theorem poses severe limitations on possible dynamics of the Hamiltonian systems. For instance, it precludes the existence of dynamic attractors. Thus, there exists a crucial qualitative difference between Hamiltonian and non-Hamiltonian dynamics. This paper considers the single-degree-of-freedom system with periodic external forcing. If such model is Hamiltonian, its stroboscopic map is expected to contain well-known generic features such as stability islands, separatrix chaos, non-attractive ``chaotic sea" etc \cite{zaslavsky2007physics}. However, the system with forcing and dissipation usually demonstrates eventual attraction of all phase trajectories to different type of attractors ($\omega$-sets) such as periodic orbits, quasiperiodic tori, or strange attractors, chaotic or non-chaotic \cite{lorenz, guckenheimer2013nonlinear, strogatz2018nonlinear}.

The force of dry friction is strictly dissipative, if the velocity is nonzero. At the same time, we are going to demonstrate that the forced vibro-impact oscillator with dry friction exhibits mixed dynamics. Namely, for the same sets of parameters, some open sets of initial conditions lead to non-attractive dynamics, similar to that of globally Hamiltonian systems. If the initial conditions lie beyond these sets, the phase trajectories are attracted to well-defined $\omega$-sets. It will be also demonstrated that both factors (impacts and friction) are crucial for appearance of the mixed dynamics.

We note that the phemenon of mixed dynamics is not a novel discovery. Similar behaviors have been observed, for example, in time-reversible systems\cite{politi,  sprott} and Fermi-Ulam model with drag\cite{Leonel_2006}. However, in the current paper we present a class of simple mechanical models where, to the best of the authors' knowledge, the phenomenon of mixed global dynamics has not been discovered before. 

The paper is as organized follows. In Section~\ref{sec:mixed}, the phenomenon of mixed dynamics is presented and illustrated by numerical simulations. Section~\ref{sec:analysis} contains analysis of the observed phenomena. Possible generalizations of the models with mixed dynamics are discussed in Section~\ref{sec:extensions}, followed by concluding remarks.

\section{\label{sec:mixed}Mixed Dynamics of Forced Vibro-impact Oscillator}
We consider a single-degree-of-freedom unit mass particle placed between two rigid walls at $x = r$ and $x = l$ with $r - l = R > 0$. Between the walls the particle moves under the action of external sinusoidal force  of period $T=2\pi/\omega$ and the Amonton-Coulomb dry friction force:
\begin{equation}\label{main}
\ddot{x} + f \; \text{sgn}\left(\dot x\right) = F \cos{(\omega t)}, \qquad  l < x < r. % - f_{\text{dry}}(\dot{x}).
\end{equation}
Strictly speaking the solution of equation~\eqref{main} should be understood in Filippov's sense, {\it i.e.} the nonlinearity is the convexification of the discontinuous sign function. We assume that the amplitude of external forcing exceeds the value of the kinetic friction, $F > f$, since otherwise any motion comes to stop.
%The dry friction term $f_{dry}(\dot{x})$ is given by
%\begin{equation}
%f_{\text{dry}}(\dot{x}) = \begin{cases} \mu, & \dot{x}=0 \\
% f \; \text{sgn}(\dot{x}), & \dot{x} \neq 0,
%\end{cases}
%\end{equation}
%where $0 < f < F$ and $\mu$ takes values between $-f$ and $f$ which are necessary for equilibrium.
We also assume the ideal reflection off the walls, {\it i.e.} the unit restitution coefficient. In other words, when the particle hits either wall at a moment $t$ the instantaneous transformation of velocity is given by
\begin{equation}\label{impact}
\dot{x}(t^-) = - \dot{x}(t^+)\quad \text{ for }\quad x(t) = r, l.
\end{equation}

Dynamics of system~\eqref{main}--\eqref{impact}, and its dependence on parameters, will be illustrated by the global phase portrait of the stroboscopic time $T$ map $\Phi$ that takes a point $\left(x(0),\; \dot x(0)\right)$ to the point $\left(x\left(T\right),\; \dot x\left(T\right)\right)$.
In particular, fixed points of $\Phi$ correspond to periodic solutions with the frequency of external forcing $\omega$. We note that rescaling of length and time can be used to normalize two out of the four system parameters, $F$, $f$, $\omega$, $R$, to $1$.

\begin{figure}[H]
\begin{subfigure}[t]{0.44\textwidth}
\includegraphics*[width=\textwidth]{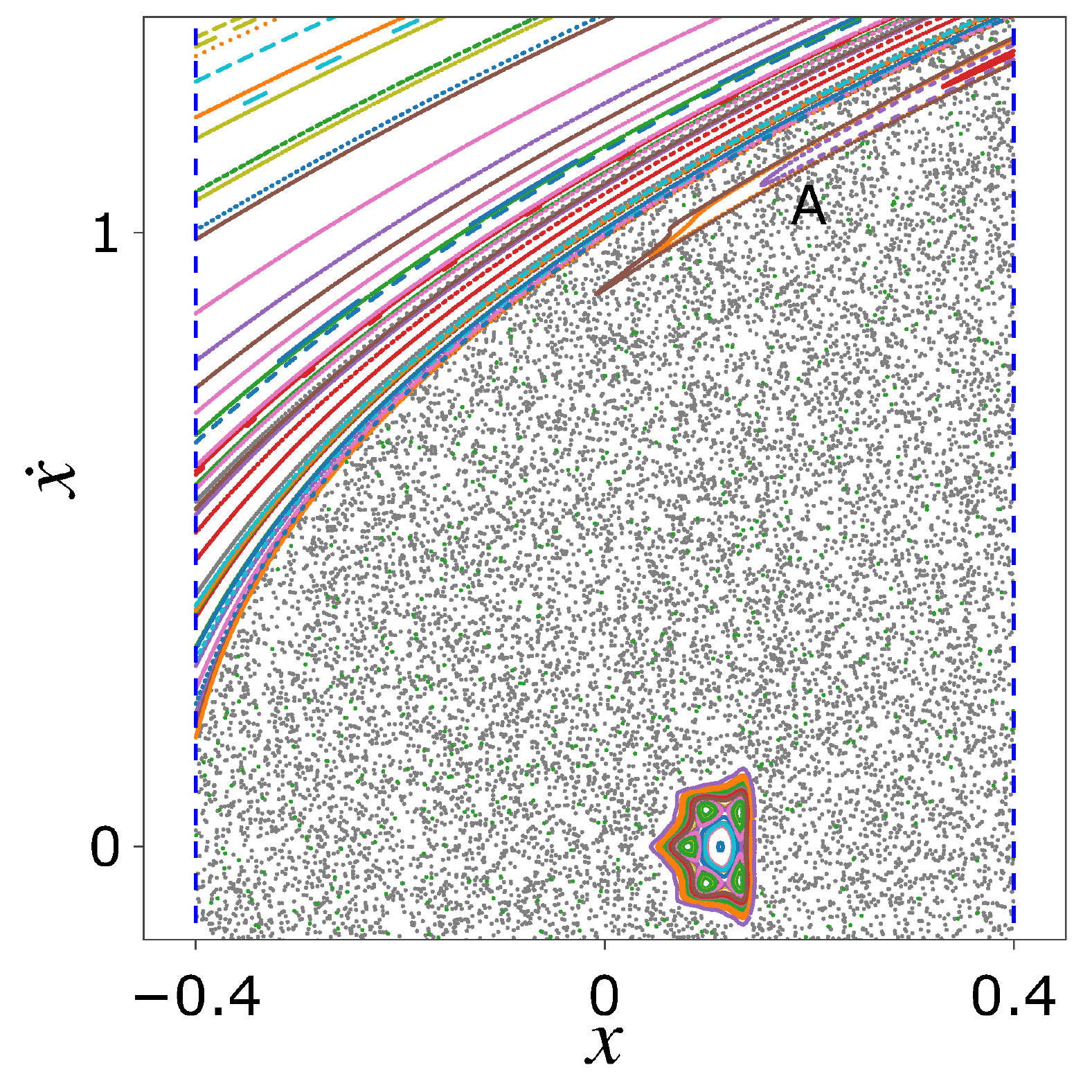}
\caption{\label{fig:narrow}}
\end{subfigure}\hfill
\begin{subfigure}[t]{0.44\textwidth}
\includegraphics*[width=\textwidth]{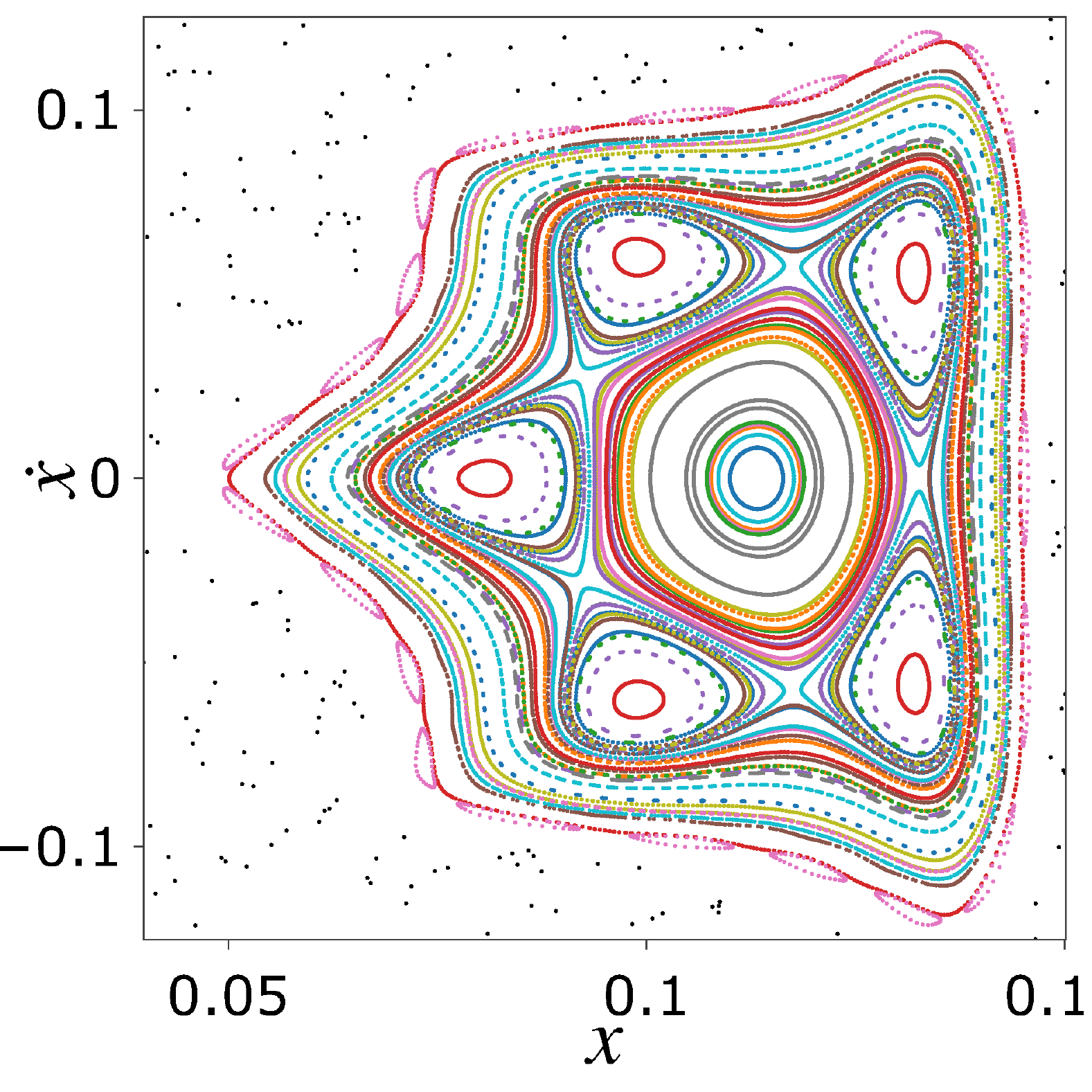}
\caption{\label{fig:zoom}}
\end{subfigure}
\begin{subfigure}[t]{0.44\textwidth}
\includegraphics*[width=\textwidth]{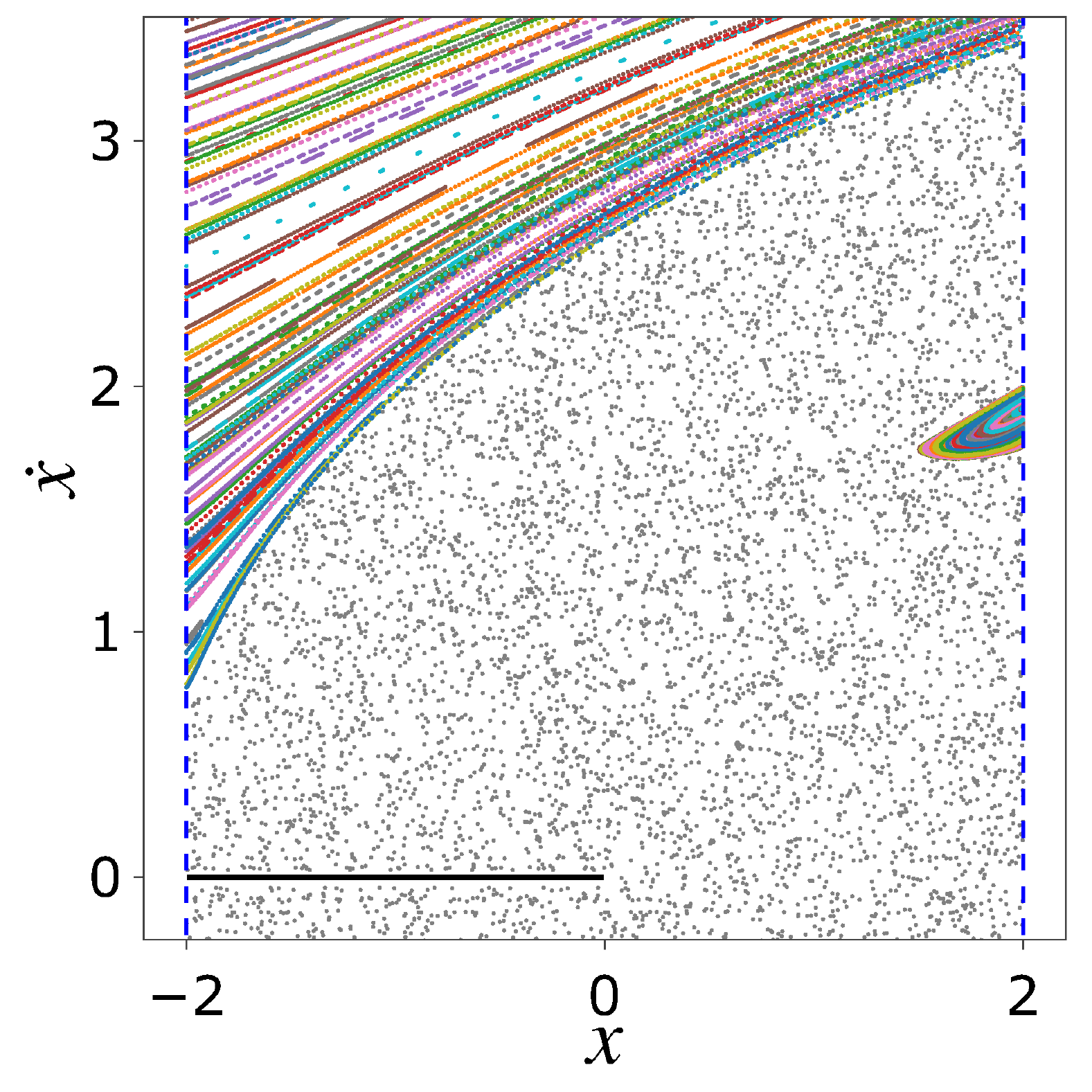}
\caption{\label{fig:wide:a}}
\end{subfigure}\hfill
\begin{subfigure}[t]{0.44\textwidth}
\includegraphics*[width=\textwidth]{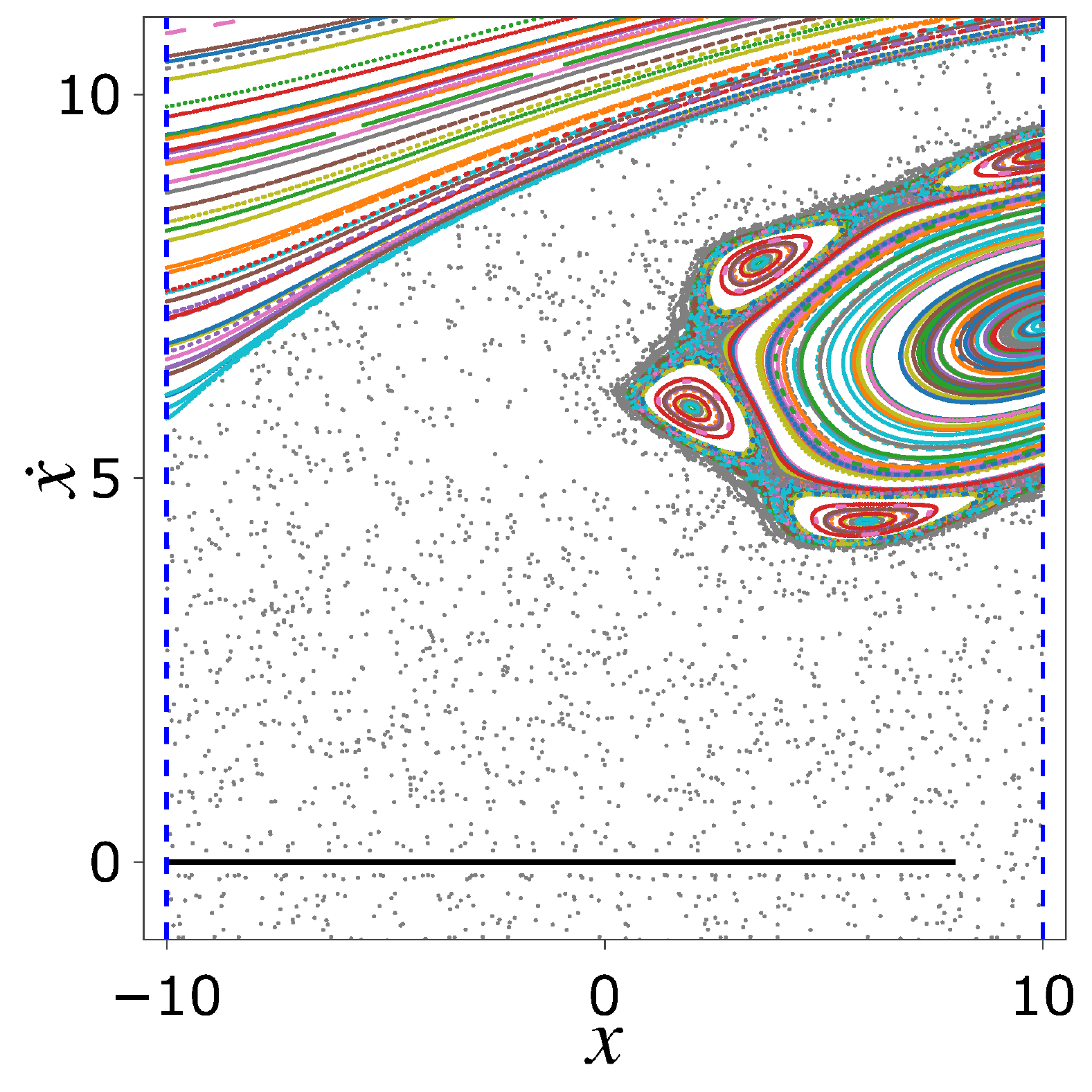}
\caption{\label{fig:wide:b}}
\end{subfigure}
\caption{\small Phase portraits of the time $T$ map $\Phi$ for $F=1$, $f=0$, $\omega = 1$ and different values of the distance between the walls, $R$. Panel (a): $R = 0.8$. Panel (b): zoom of an invariant island from Panel (a). Panel (c): $R = 4$. Panel (d): $R = 20$. Blue dashed lines denote the walls; colored lines represent invariant curves (tori); gray dots  correspond to the chaotic sea; the black line segment consists of fixed points corresponding to periodic solutions without impacts. \label{fig:nofriction}}
\end{figure}

Figure~\ref{fig:narrow} presents the phase portrait of the stroboscopic map $\Phi$ of system~\eqref{main}--\eqref{impact} in the frictionless case, $f = 0$ (other system parameters are defined in the figure caption). This Hamiltonian system can be used as a reference for further results. An invariant island in the center contains a fixed point and a 5-periodic orbit of the map $\Phi$, which are surrounded by invariant curves representing quasiperiodic motions, see Figure~\ref{fig:zoom}.
This fixed point corresponds to a periodic solution of the fundamental frequency $\omega$ with four impacts at the walls and two additional turning points per period. % and a subharmonic periodic orbit of frequency $5\omega$
Another region of quasiperiodic motions is observed for large velocities. Between these two regions there is a wide domain of chaotic dynamics, see Figure~\ref{fig:narrow}. The third island of quasiperiodic motions,
which is indicated by letter $A$ on Figure~\ref{fig:narrow}, is ``centered" at a fixed point corresponding to a $2\pi/\omega$-periodic motion with two impacts per period and no other turning points. The latter trajectories belong to the class of {\em non-sticking} solutions, {\it i.e.} solutions which have non-zero velocity at all times. They are important for further analysis.

\begin{figure}[H]
\centering
\begin{subfigure}[t]{0.44\textwidth}
\includegraphics*[width=\textwidth]{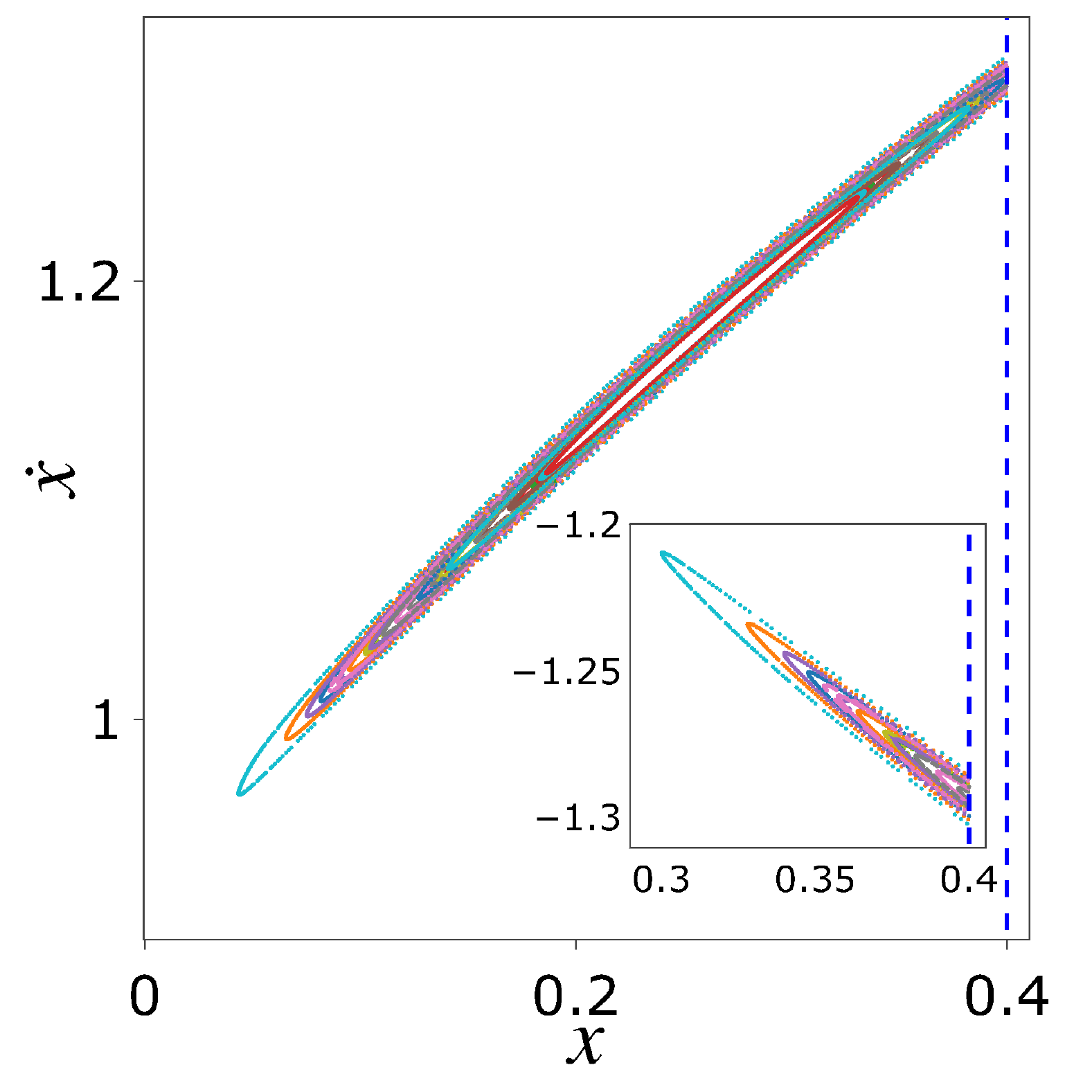}
\caption{\label{fig:narrow:0005a}}
\end{subfigure}
\hfill
\begin{subfigure}[t]{0.44\textwidth}
\includegraphics*[width=\textwidth]{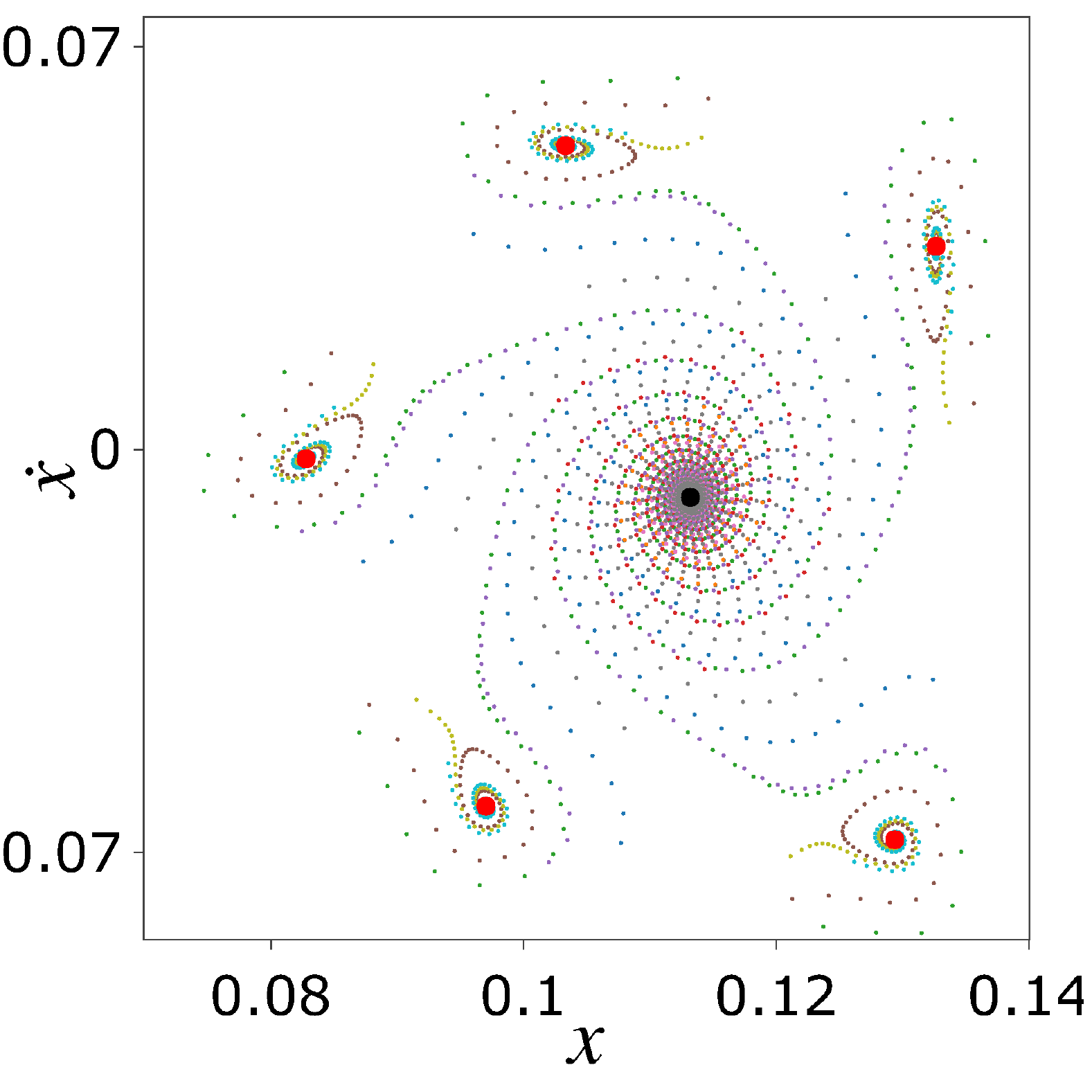}
\caption{\label{fig:zoom-with-friction}}
\end{subfigure}
\caption{\small Same as in Figures~\ref{fig:nofriction}a,b, for $f = 0.005$. All the trajectories starting from the white region are attracted to the asymptotically stable focus denoted by the black dot on Panel (b). Red dots correspond to the stable $5$-periodic orbit, which can be traced back by continuation to the neutrally stable orbit of the same period in Figure~\ref{fig:narrow-with-friction}b.\label{fig:narrow-with-friction}}
\end{figure}

\begin{figure}[H]
\begin{subfigure}[t]{0.31\textwidth}
\includegraphics[width=\textwidth]{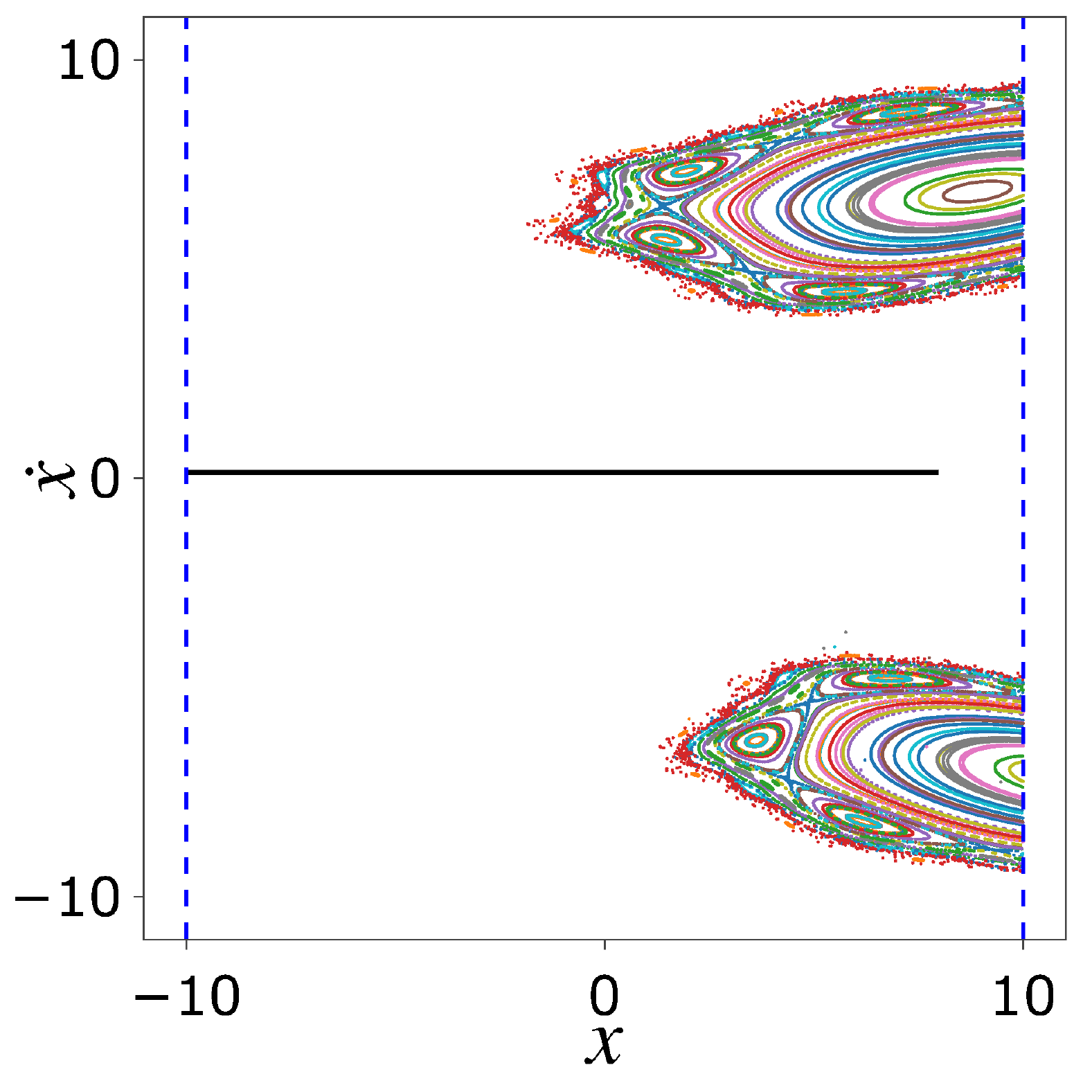}
\caption{}
\end{subfigure}
\hfill
\begin{subfigure}[t]{0.31\textwidth}
\includegraphics[width=\textwidth]{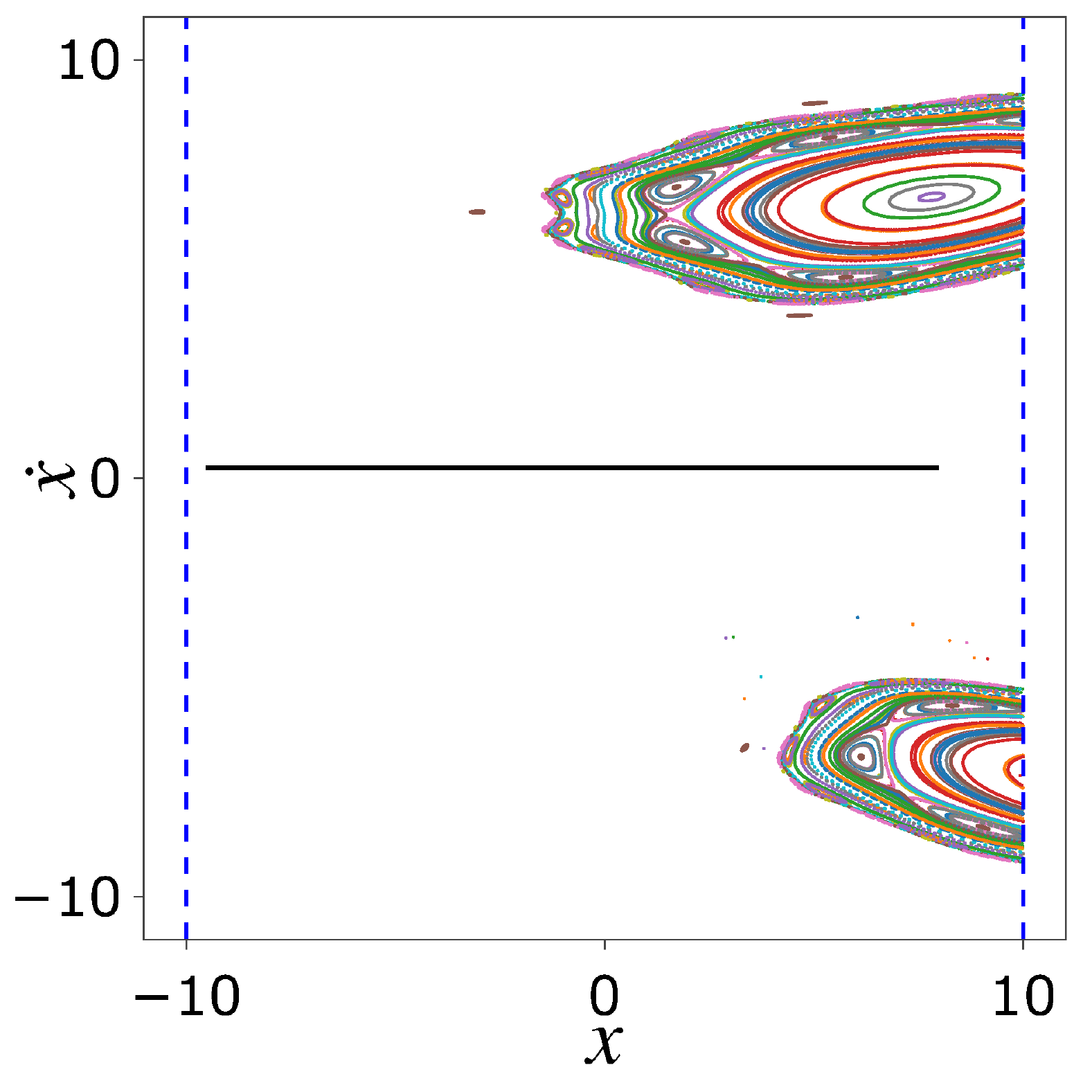}
\caption{}
\end{subfigure}
\hfill
\begin{subfigure}[t]{0.31\textwidth}
\includegraphics[width=\textwidth]{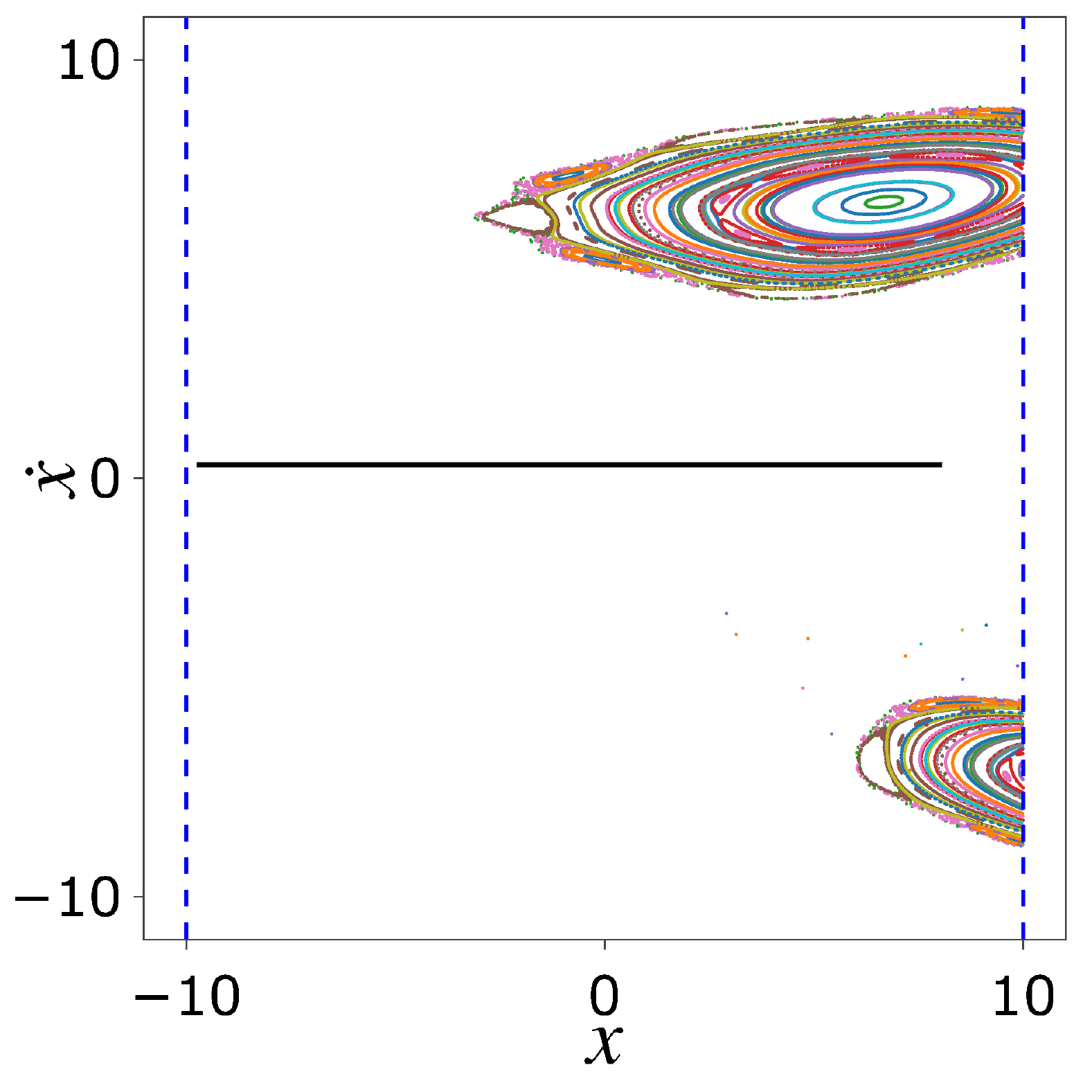}
\caption{}
\end{subfigure}
\begin{subfigure}[t]{0.31\textwidth}
\includegraphics[width=\textwidth]{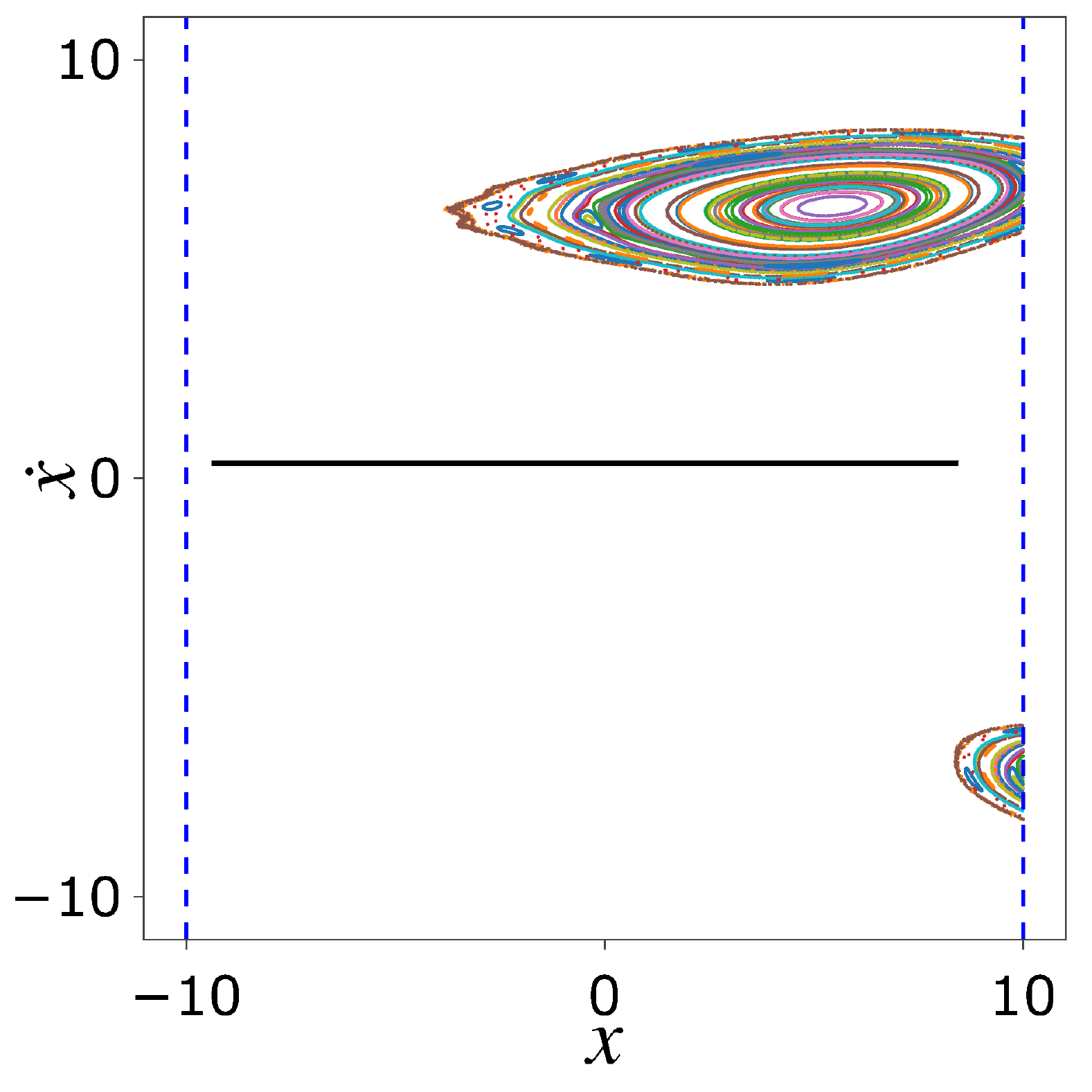}
\caption{}
\end{subfigure}
\hfill
\begin{subfigure}[t]{0.31\textwidth}
\includegraphics[width=\textwidth]{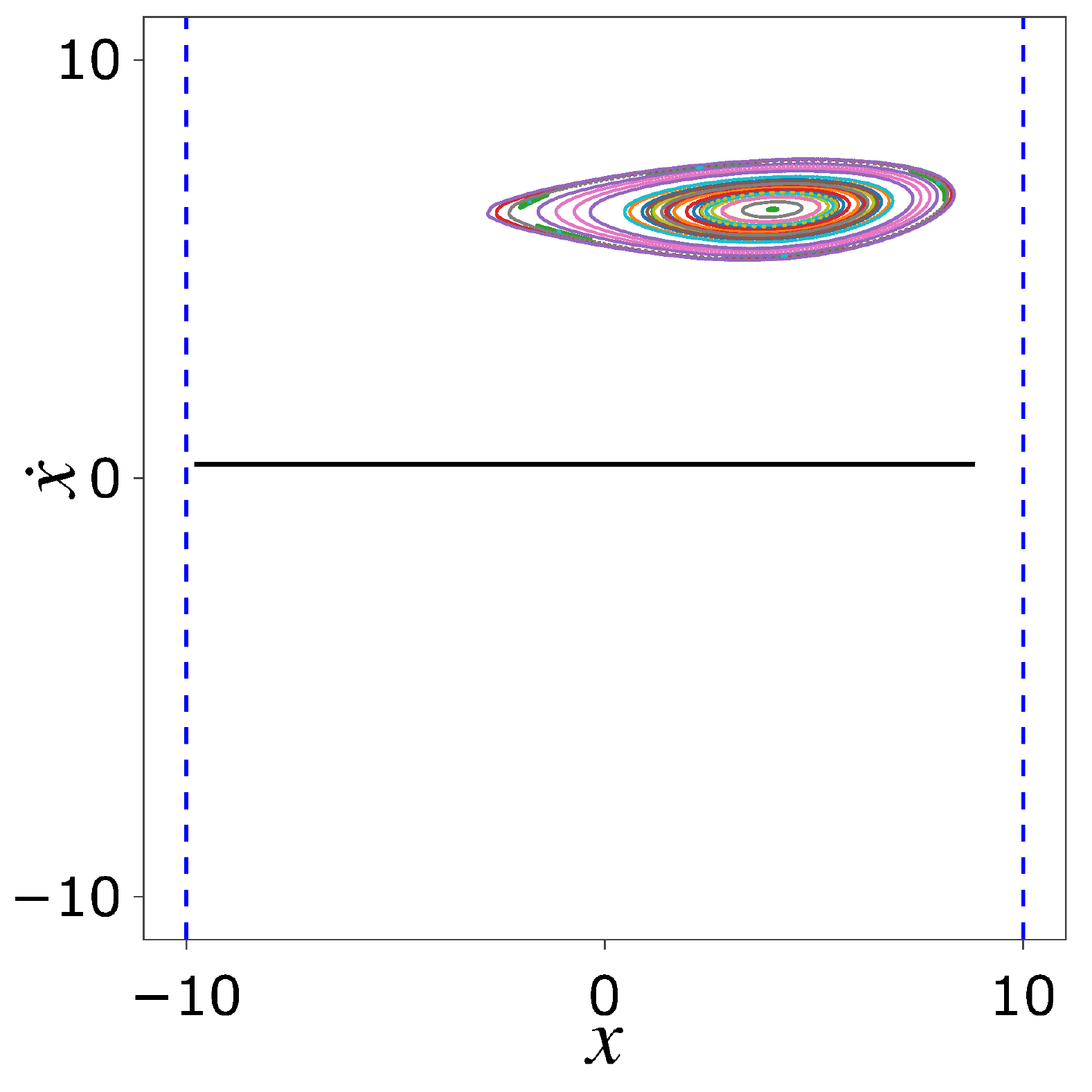}
\caption{}
\end{subfigure}
\hfill
\begin{subfigure}[t]{0.31\textwidth}
\includegraphics[width=\textwidth]{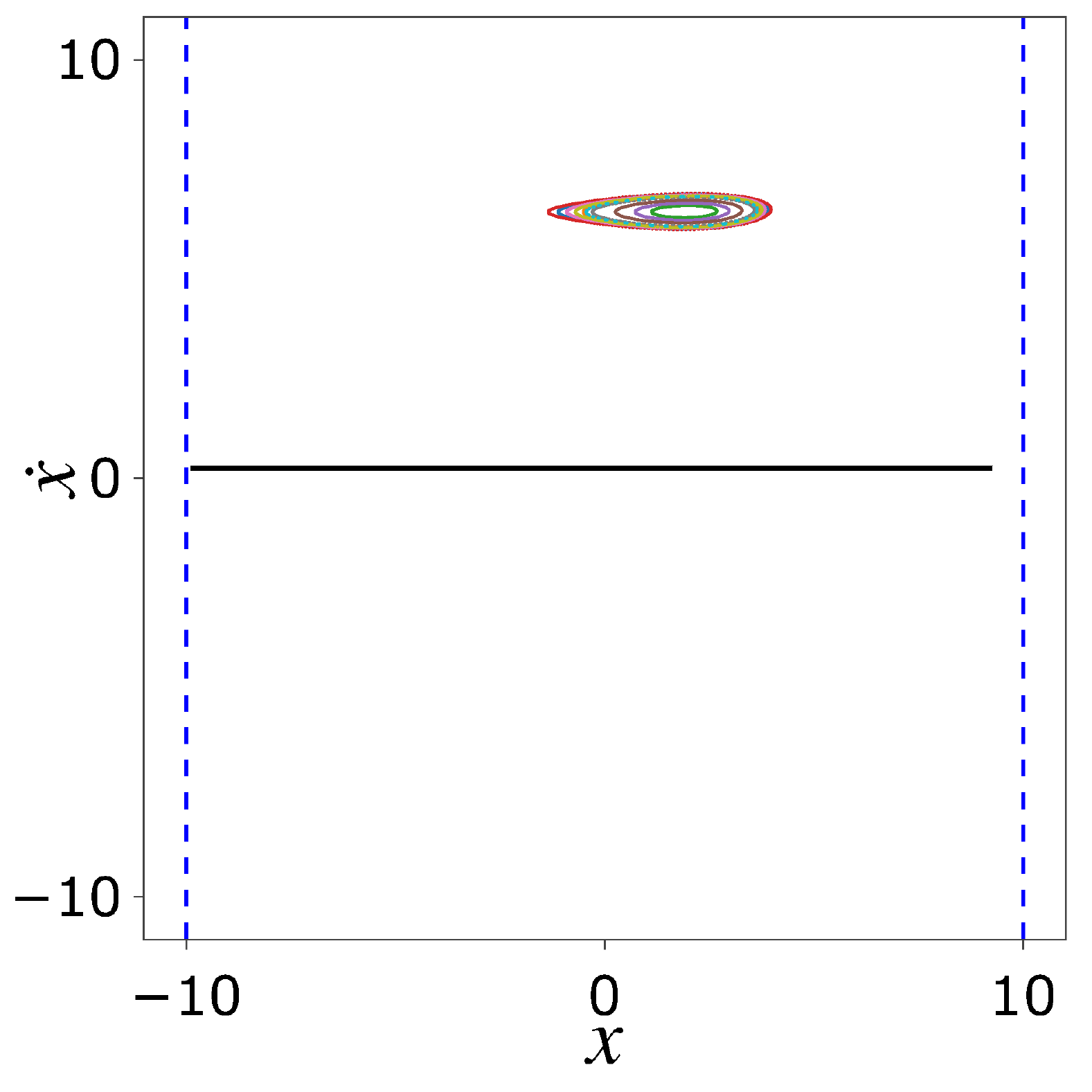}
\caption{}
\end{subfigure}
\caption{\small Evolution of the phase portrait with increasing friction. Panels (a)--(f) correspond to the values of $f$ from $0.1$ to $0.6$ with increment $0.1$, respectively. Other parameters are the same as in Figure~\ref{fig:wide:b}.\label{fig:many:wide}}
\end{figure}

The above phase portrait is slightly transformed for larger values of the parameter $R = r-l$. In particular, when $R > F/\omega^2$, the system has periodic solutions without impacts,
\begin{equation}\label{noimpacts}
x(t) = -F{\omega^{-2}} \cos {(\omega t)} + C,
\end{equation}
which are represented by a horizontal segment of fixed points in Figures~\ref{fig:nofriction}c,d.
Introducing small friction $f = 0.005$ with the parameter set used in Figure~\ref{fig:narrow} has a significant impact on the phase portrait. The invariant region of Hamiltonian dynamics which is marked by $A$ on Figure~\ref{fig:narrow} persists after this perturbation, see Figure~\ref{fig:narrow-with-friction}. As shown below the time $T$ map is area-preserving within this invariant region. However, the other regions of quasi-periodic dynamics and chaos of Figure~\ref{fig:narrow} are now all replaced by the basins of attraction of the fixed point and period $5$ orbit of the map $\Phi$ on Figure~\ref{fig:narrow-with-friction}. These asymptotically stable fixed point and period $5$ orbit can be traced back by continuation to the neutrally stable fixed point and period 5 trajectory on Figure~\ref{fig:zoom}. In other words, in this case we observe two periodic attractors and, simultaneously, an invariant region within which the area is preserved.
Figure~\ref{fig:many:wide} presents the transformation  of the phase portrait shown in Figure~\ref{fig:wide:b} with increasing friction. Here again the invariant island of quasi-periodic dynamics ``centered" at the fixed point, which corresponds to a non-sticking periodic solution, persists after the introduction of friction. This island coexists with an attractor which in this case consists of the periodic solutions without impacts that can be traced back to solutions~\eqref{noimpacts}. This attractor is represented by the black horizontal line segment on Figure~\ref{fig:many:wide}. In particular, its basin of attraction replaces the chaotic regions and the invariant tori with large velocities shown in Figure~\ref{fig:wide:b}, even in the case of small friction.

\begin{figure}[H]
\centering
\includegraphics*[width=0.4\textwidth]{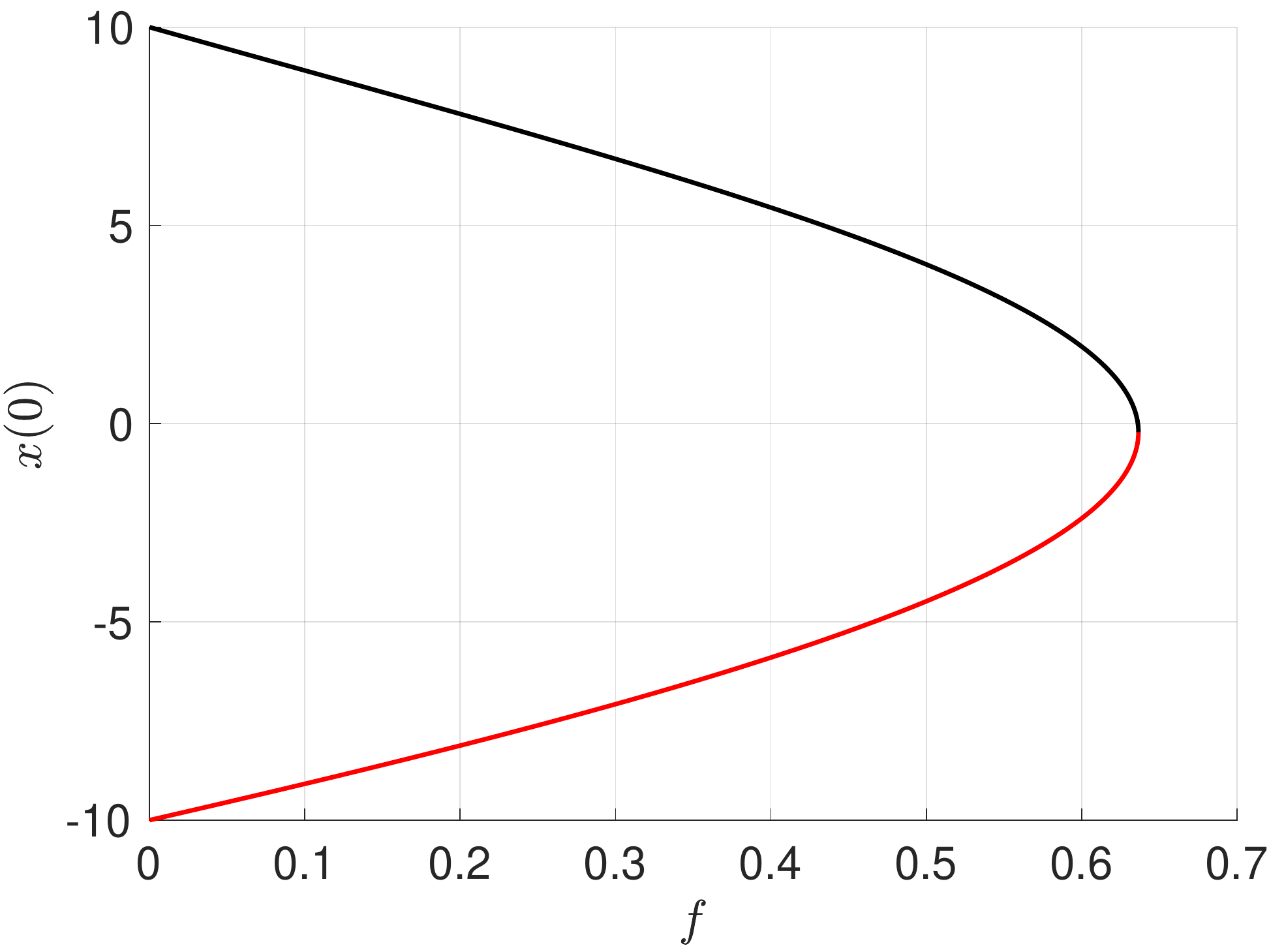}
\caption{\small Saddle-center bifurcation of non-sticking periodic orbits. The black part corresponds to neutrally stable periodic solutions (centers), the red part corresponds to periodic solutions of saddle type. All the parameters except $f$ are the same as in Figure~\ref{fig:wide:b}. This bifurcation diagram was obtained using \sc{coco}\cite{danko}.\label{fig:fold}}
\end{figure}

\begin{figure}[H]
\centering
\begin{subfigure}[t]{0.48\textwidth}
\includegraphics*[width=\textwidth]{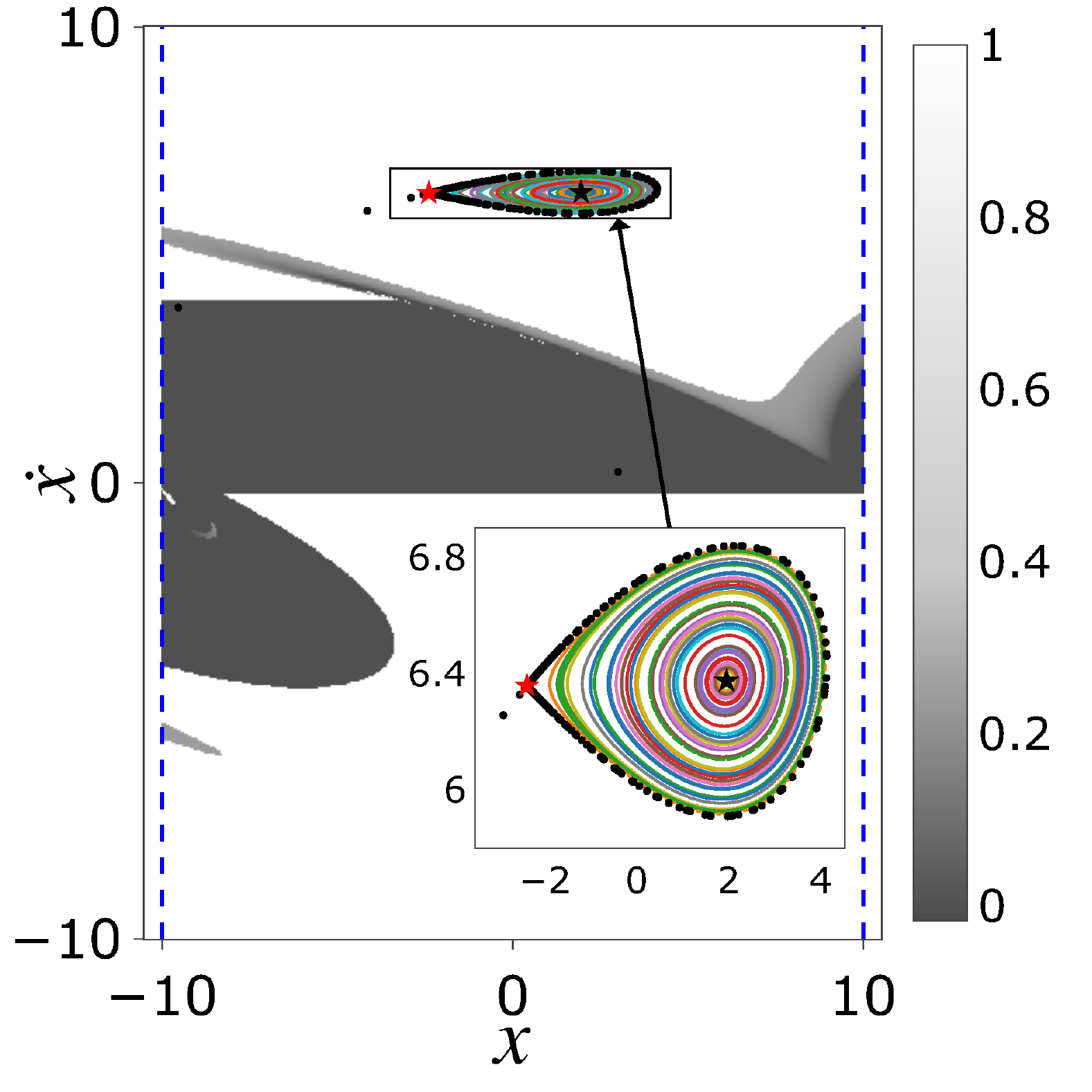}
\caption{\label{fig:regions}}
\end{subfigure}
\hfill
\begin{subfigure}[t]{0.48\textwidth}
\includegraphics*[width=\textwidth]{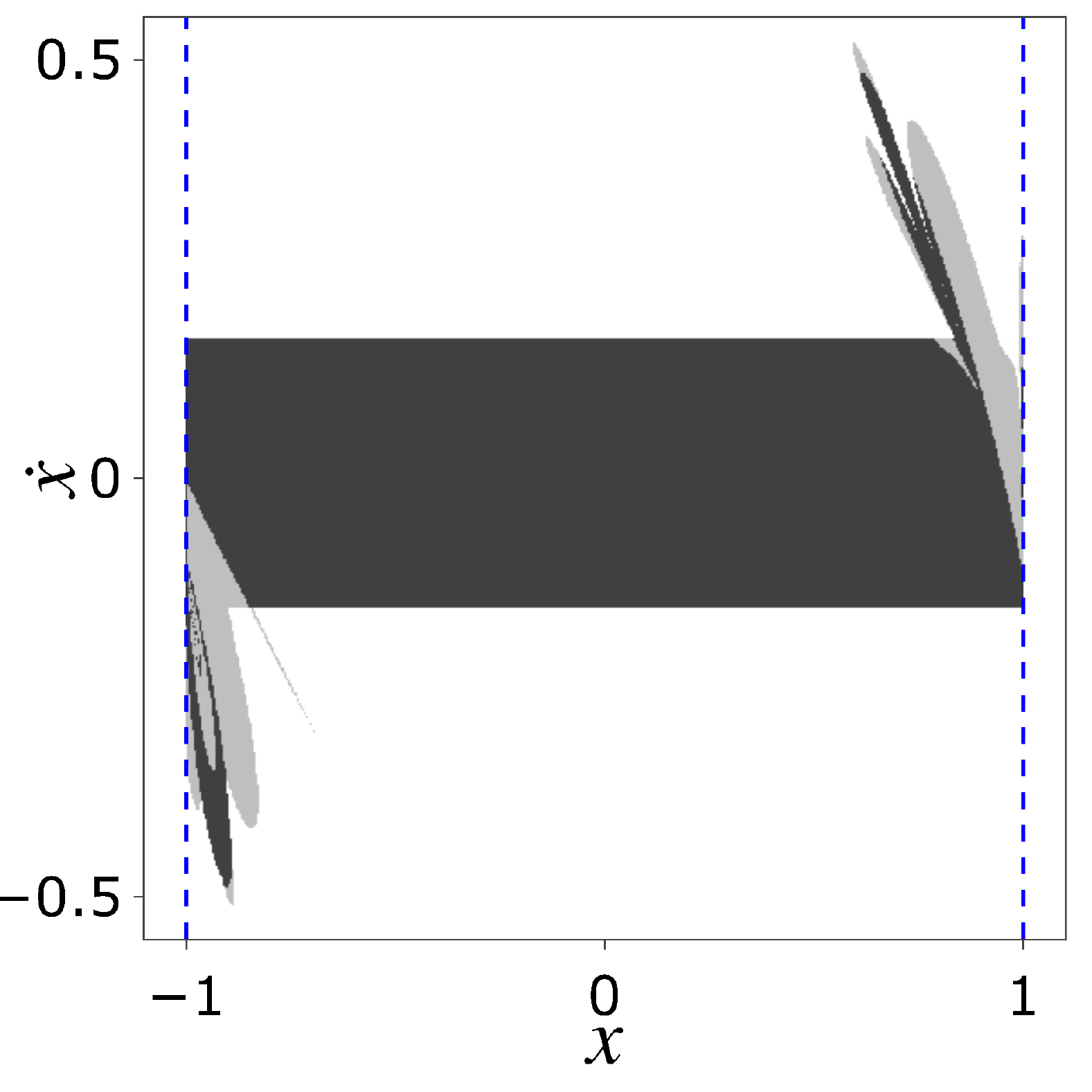}
\caption{\label{fig:regs}}
\end{subfigure}
\caption{\small Panel (a): White and gray regions correspond to parts of $\Omega_{pres}$ and $\Omega_{dissip}$, respectively. The shades of gray in $\Omega_{dissip}$ are scaled according to the values of $\det{(\Phi^\prime)}$. Red star corresponds to the unstable fixed point of $\Phi$. This orbit has eigenvalues $\lambda_1 = 0.3159$ and $\lambda_2 = 3.1659$. A perturbation of the initial condition along the unstable direction produces a trajectory (black dots) which follows the boundary of the invariant island. This trajectory converges to a fixed point corresponding to a periodic solution without impacts. Black star corresponds to the stable fixed point. Panel (b): White and gray regions correspond to $\Omega_{pres}$ and $\Omega_{dissip}$, respectively. The image $\Phi(\Omega_{dissip})$ is dark gray. We can see that $\Phi(\Omega_{dissip}) \subseteq \Omega_{dissip}$.\label{fig:fig}}
\end{figure}

Figure~\ref{fig:many:wide} shows that the invariant island of Hamiltonian dynamics shrinks with increasing kinetic friction $f$ until it disappears in a saddle-center bifurcation of two non-sticking periodic orbits  at a critical kinetic friction value $f = f_{\text{crit}}$, see Figure~\ref{fig:fold}. For $f > f_{\text{crit}}$ the set of periodic solutions without impacts becomes a global attractor. For $f < f_{\text{crit}}$ near the bifurcation point, the homoclinic tangle of the fixed point, which is associated with the unstable non-sticking periodic trajectory, creates the chaotic boundary between the invariant island of Hamiltonian dynamics and the basin of attraction of periodic solutions without impacts.
The island of Hamiltonian dynamics is ``centered'' at the fixed point associated with a neutrally stable non-sticking periodic solution, see Figures~\ref{fig:three:per}a and~\ref{fig:regions}.
\begin{figure}[H]
\begin{subfigure}[t]{0.31\textwidth}
\includegraphics[width=\textwidth]{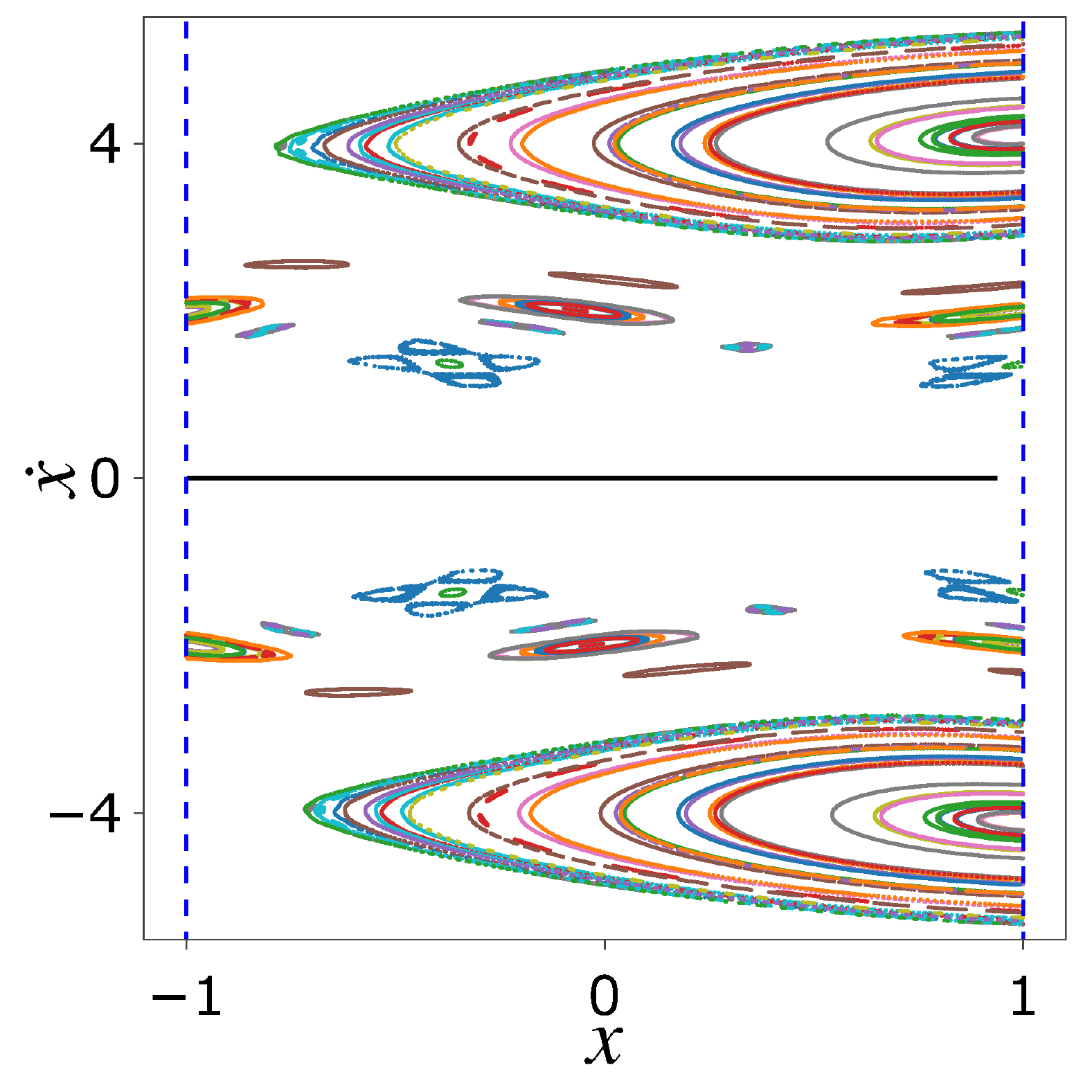}
\caption{}
\end{subfigure}
\hfill
\begin{subfigure}[t]{0.31\textwidth}
\includegraphics[width=\textwidth]{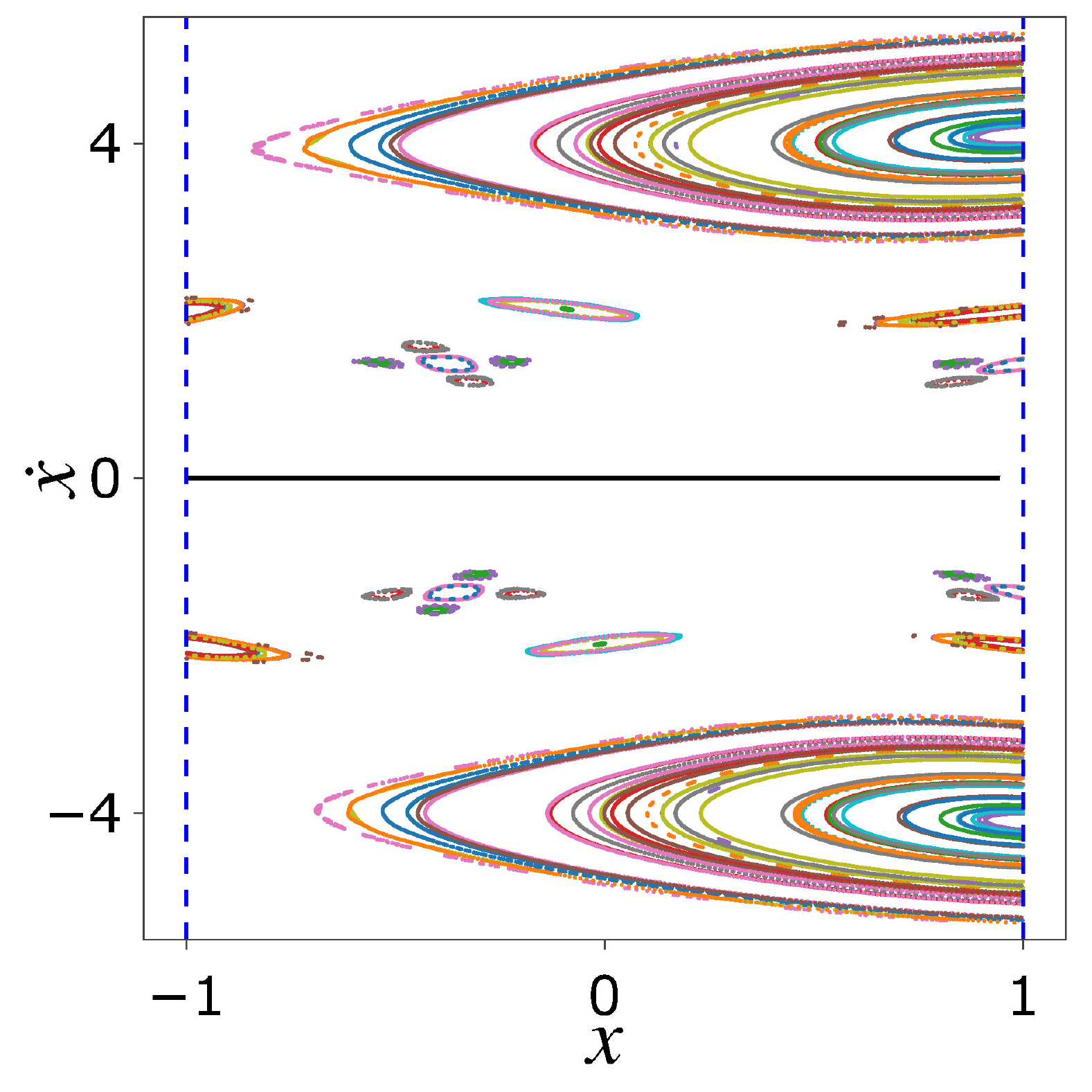}
\caption{}
\end{subfigure}
\hfill
\begin{subfigure}[t]{0.31\textwidth}
\includegraphics[width=\textwidth]{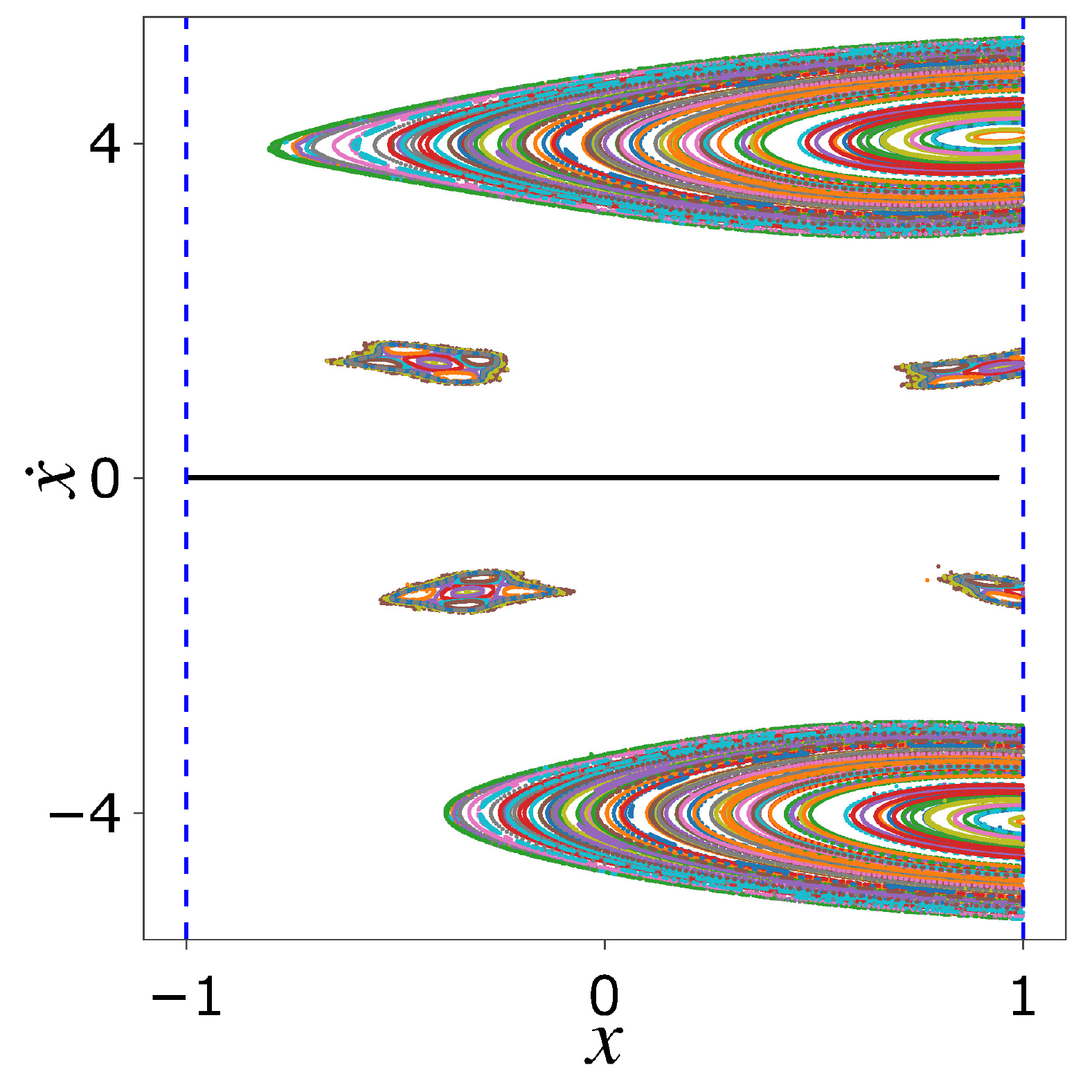}
\caption{\label{fig:f005}}
\end{subfigure}
\begin{subfigure}[t]{0.31\textwidth}
\includegraphics[width=\textwidth]{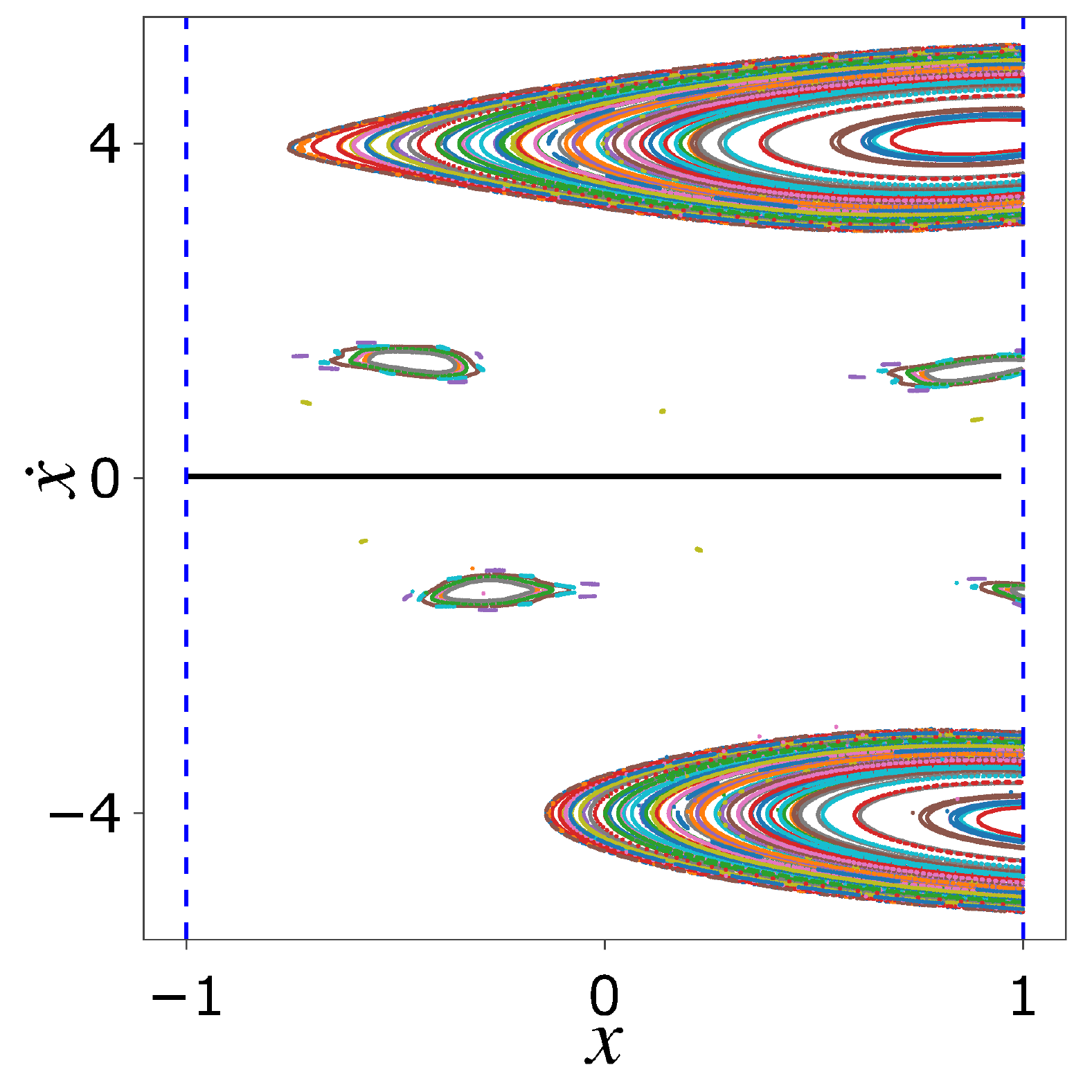}
\caption{}
\end{subfigure}
\hfill
\begin{subfigure}[t]{0.31\textwidth}
\includegraphics[width=\textwidth]{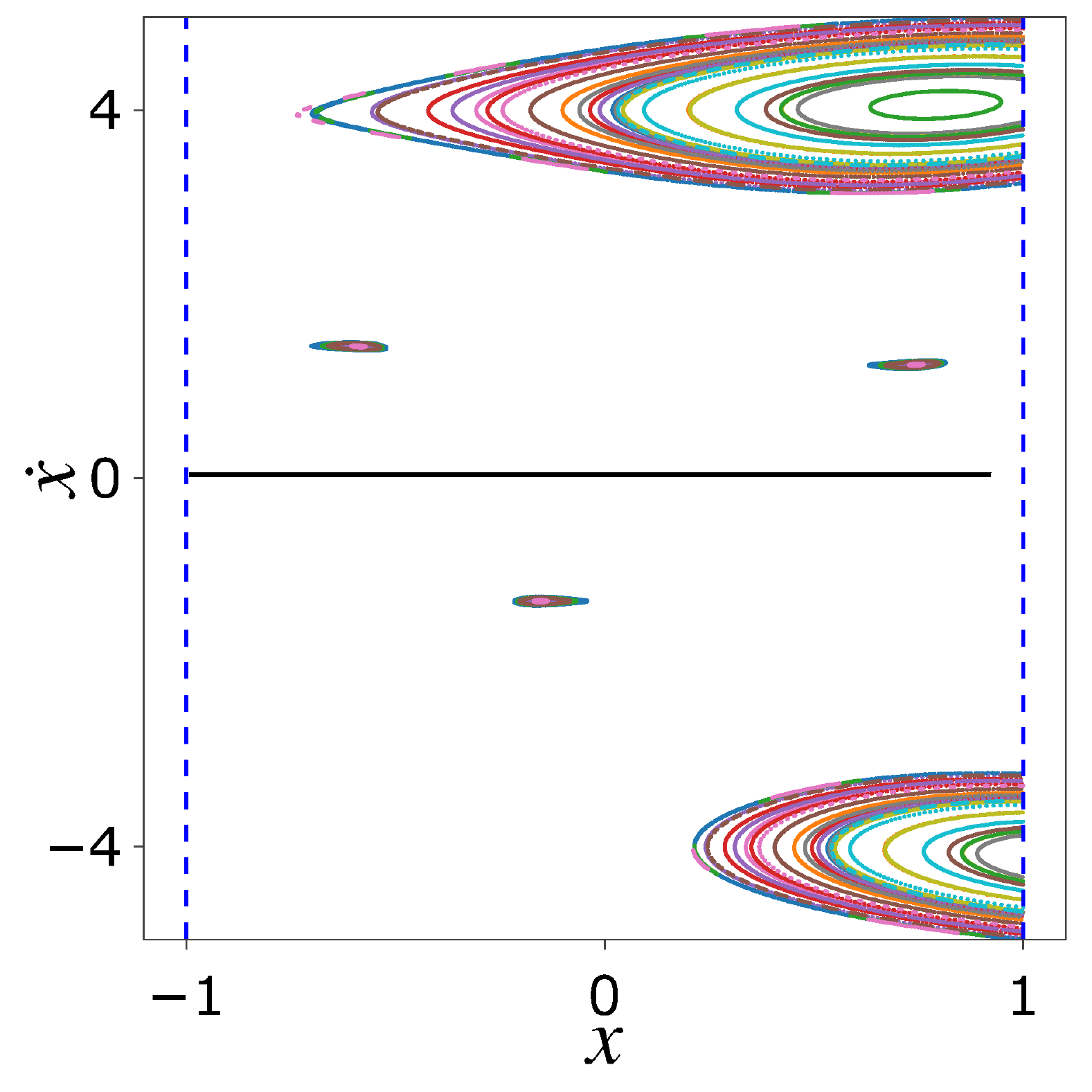}
\caption{}
\end{subfigure}
\hfill
\begin{subfigure}[t]{0.31\textwidth}
\includegraphics[width=\textwidth]{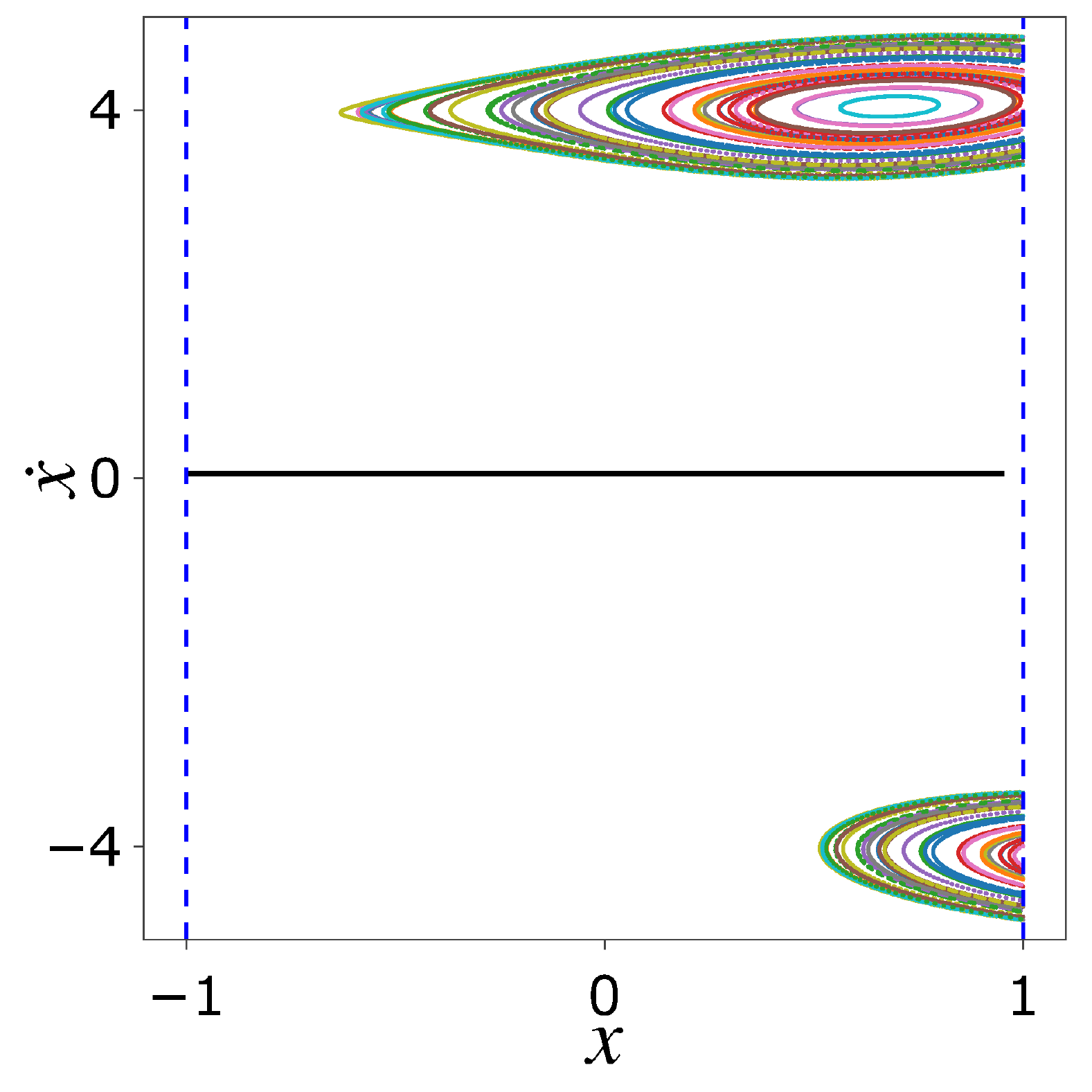}
\caption{}
\end{subfigure}
\caption{\small Same as Figure~\ref{fig:many:wide}, for $F = 1$, $\omega = 2\pi$, $l = -1$ and $r = 1$. Panels (a)--(f) correspond to the following values of $f$: $0.005$, $0.01$, $0.05$, $0.1$, $0.2$, $0.3$.\label{fig:many:2pi}}
\end{figure}
Figure~\ref{fig:many:2pi} demonstrates a similar scenario (for a different parameter set) but with multiple invariant islands of Hamiltonian dynamics. All these islands shrink with increasing friction and disappear in saddle-center bifurcations of non-sticking periodic solutions. For example, of many isolated islands shown in Figures~\ref{fig:many:2pi}a,b, only the island surrounding a non-sticking periodic solution of fundamental frequency $\omega$ and three islands surrounding a subharmonic non-sticking solution of frequency $\omega/3$ survive in Figures~\ref{fig:many:2pi}c-e corresponding to higher friction values. Further, the latter three islands disappear in Figure~\ref{fig:many:2pi}f corresponding to even higher friction $f$.

\begin{figure}[H]
\centering
\begin{subfigure}[t]{0.45\textwidth}
\includegraphics*[width=\textwidth]{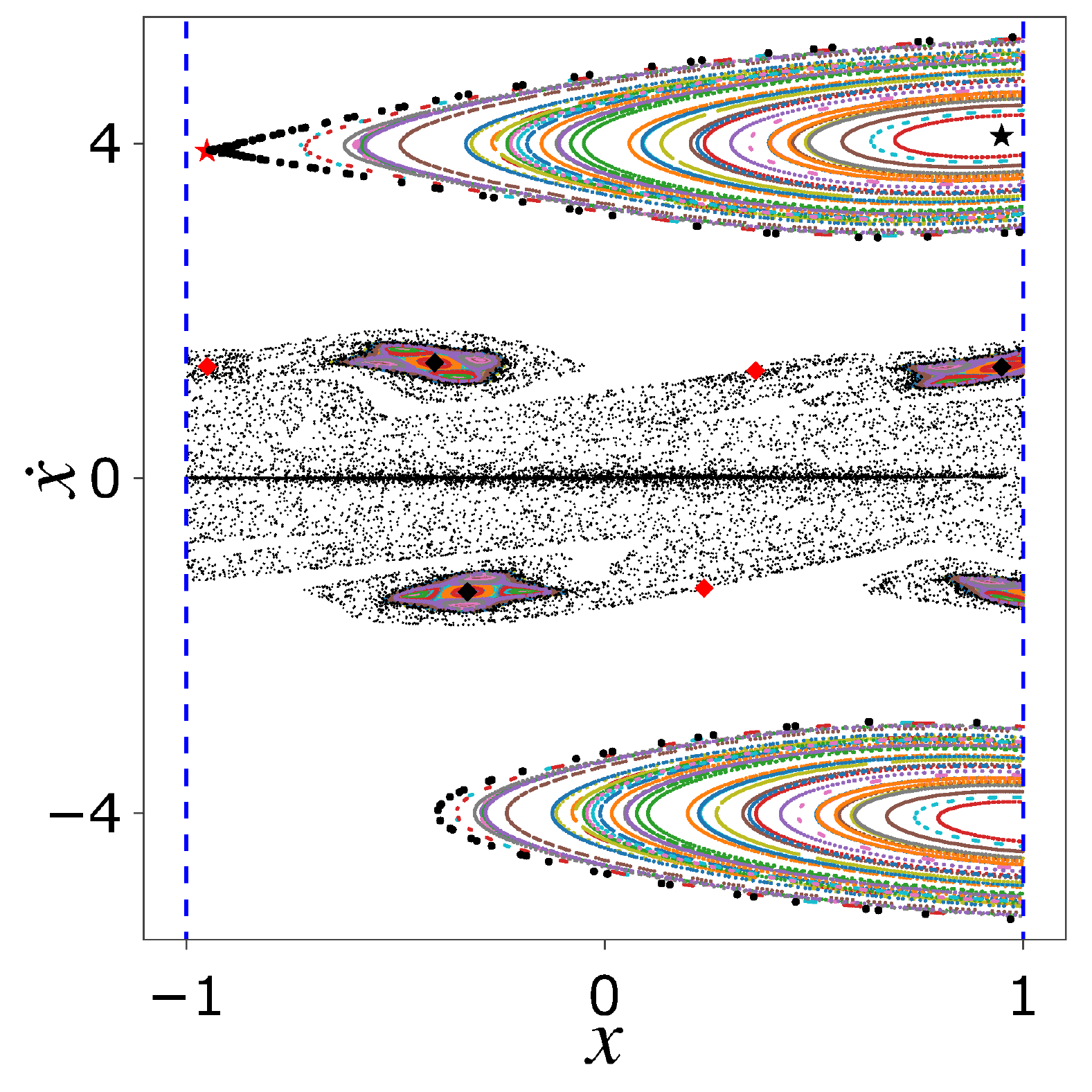}
\caption{}
\end{subfigure}
\hfill
\begin{subfigure}[t]{0.45\textwidth}
\includegraphics*[width=\textwidth]{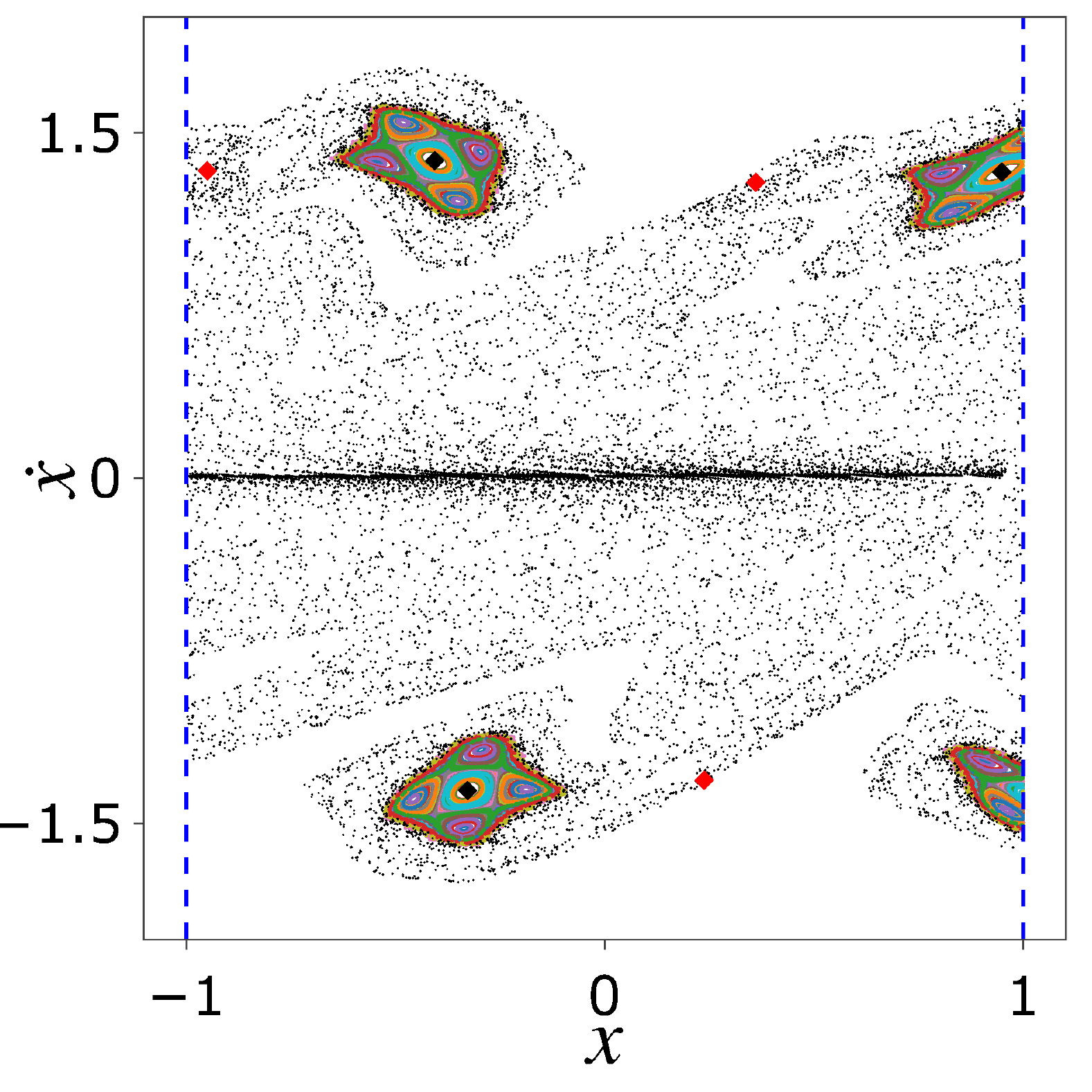}
\caption{}
\end{subfigure}
\caption{\small Black and red stars correspond to the neutrally stable (center) and unstable (saddle) fixed points of $\Phi$, respectively. Similarly, black and red diamonds are stable and unstable three-periodic orbits. Big black dots correpond to a trajectory which belongs to the homoclinic tangle bounding the invariant island of invariant curves (tori). Small black dots correspond to some trajectories with initial conditions lying in the vicinity of unstable three-periodic orbit. Such trajectories converge to the only attractor which consists of fixed points corresponding to periodic solutions without impacts. The parameters are the same as in Figure~\ref{fig:f005}.\label{fig:three:per}}
\end{figure}

\section{\label{sec:analysis}Analysis}
Hamiltonian dynamics of system~\eqref{main}--\eqref{impact} can be  explained by lifting it to the unconstrained single-degree-of-freedom Hamiltonian system
\begin{equation}\label{hamilsys}
\ddot{q} + \frac{\partial V}{\partial q} = \frac{f}{R}
\end{equation}
with the space and time periodic piecewise linear potential
\begin{equation}
V(t,\; q) = \frac{F}{R}\cos{(\omega t)}W(q), \qquad W(q) = \begin{cases}
q, & 0\le q < 1,\\
2-q, & 1 \le q < 2,
\end{cases}\qquad W(q+2) = W(q).
\end{equation}
Any non-sticking solution of~\eqref{hamilsys} is mapped to a non-sticking solution of~\eqref{main}--\eqref{impact} by the simple relationship
\begin{equation}\label{map}
x = R \cdot W(q) + l.
\end{equation}
Therefore, if all solutions
%of system~\eqref{main}--\eqref{impact}
starting from a domain $\Omega$  of the state space of system~\eqref{main}--\eqref{impact} are non-sticking and $\Omega$ is invariant for the time $T$ map $\Phi$ of this system, then $\Omega$ is a region of Hamiltonian dynamics. It is important to note that systems~\eqref{main}--\eqref{impact} and~\eqref{hamilsys} are topologically conjugate only on a part of their state spaces, and trajectories with sticking are not related by equation~\eqref{map}.

\begin{figure}[H]
\centering
\begin{subfigure}[t]{0.32\textwidth}
\includegraphics*[width=\textwidth]{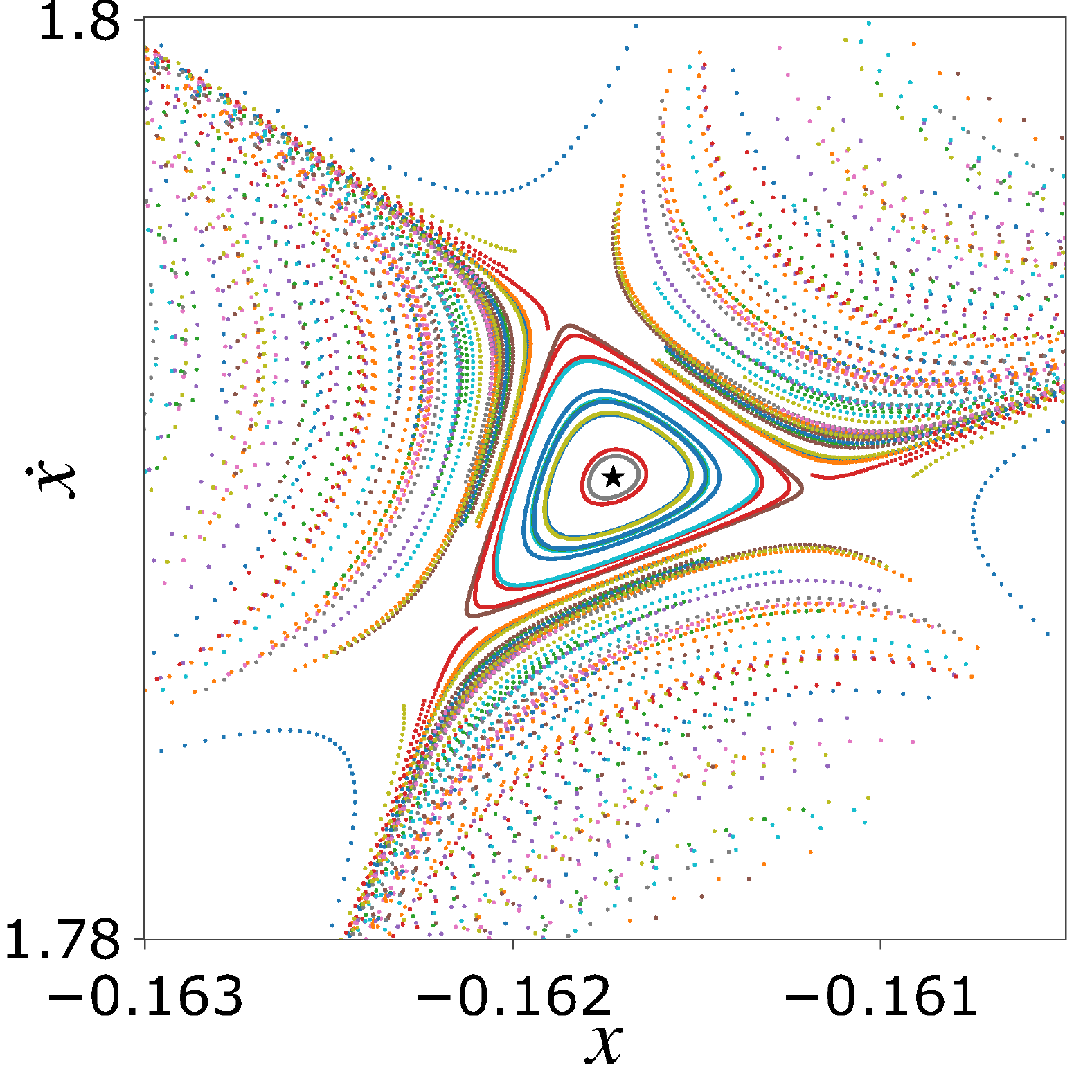}
\caption{$R = 0.89$}
\end{subfigure}
\hfill
\begin{subfigure}[t]{0.32\textwidth}
\includegraphics*[width=\textwidth]{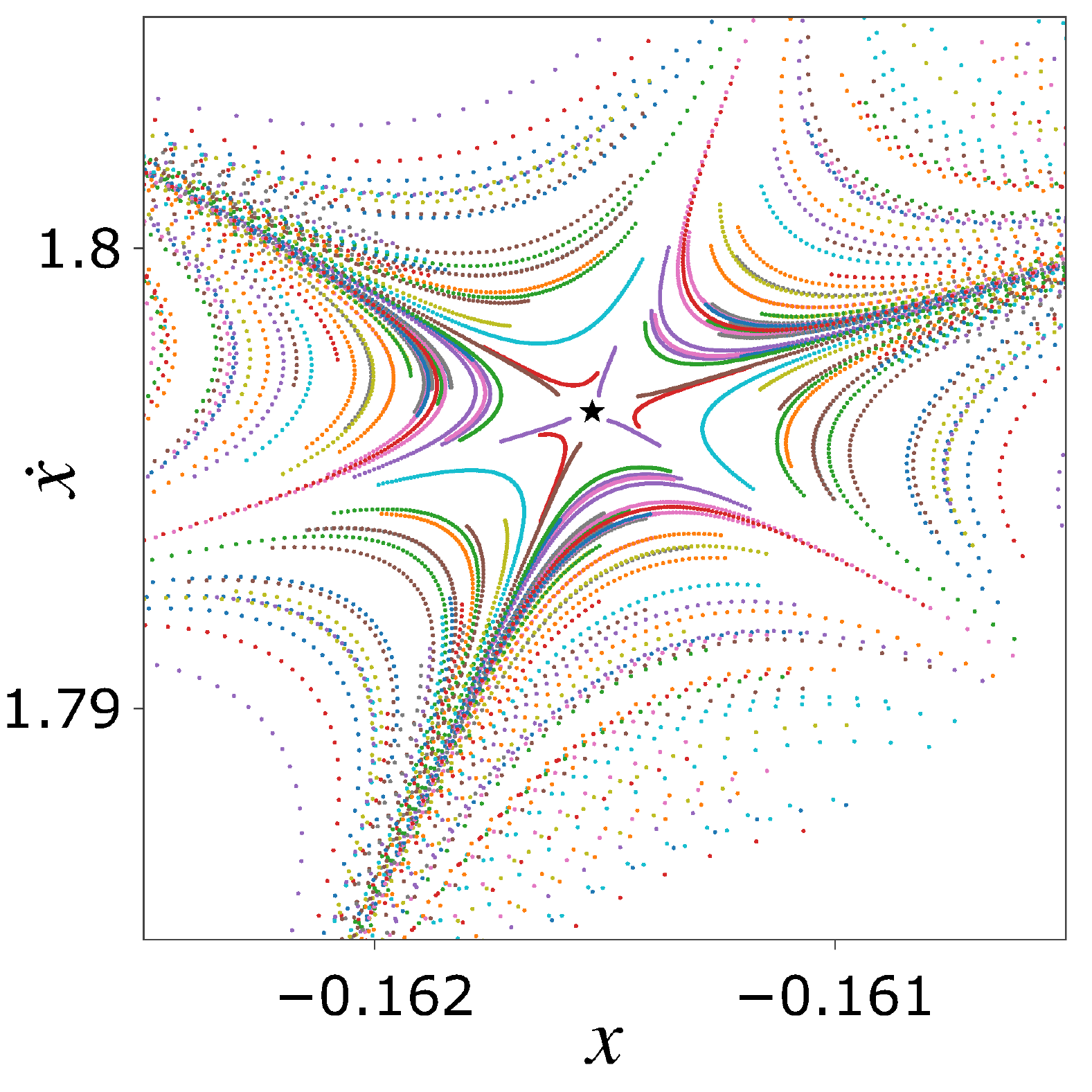}
\caption{$R = 0.8932$}
\end{subfigure}
\hfill
\begin{subfigure}[t]{0.32\textwidth}
\includegraphics*[width=\textwidth]{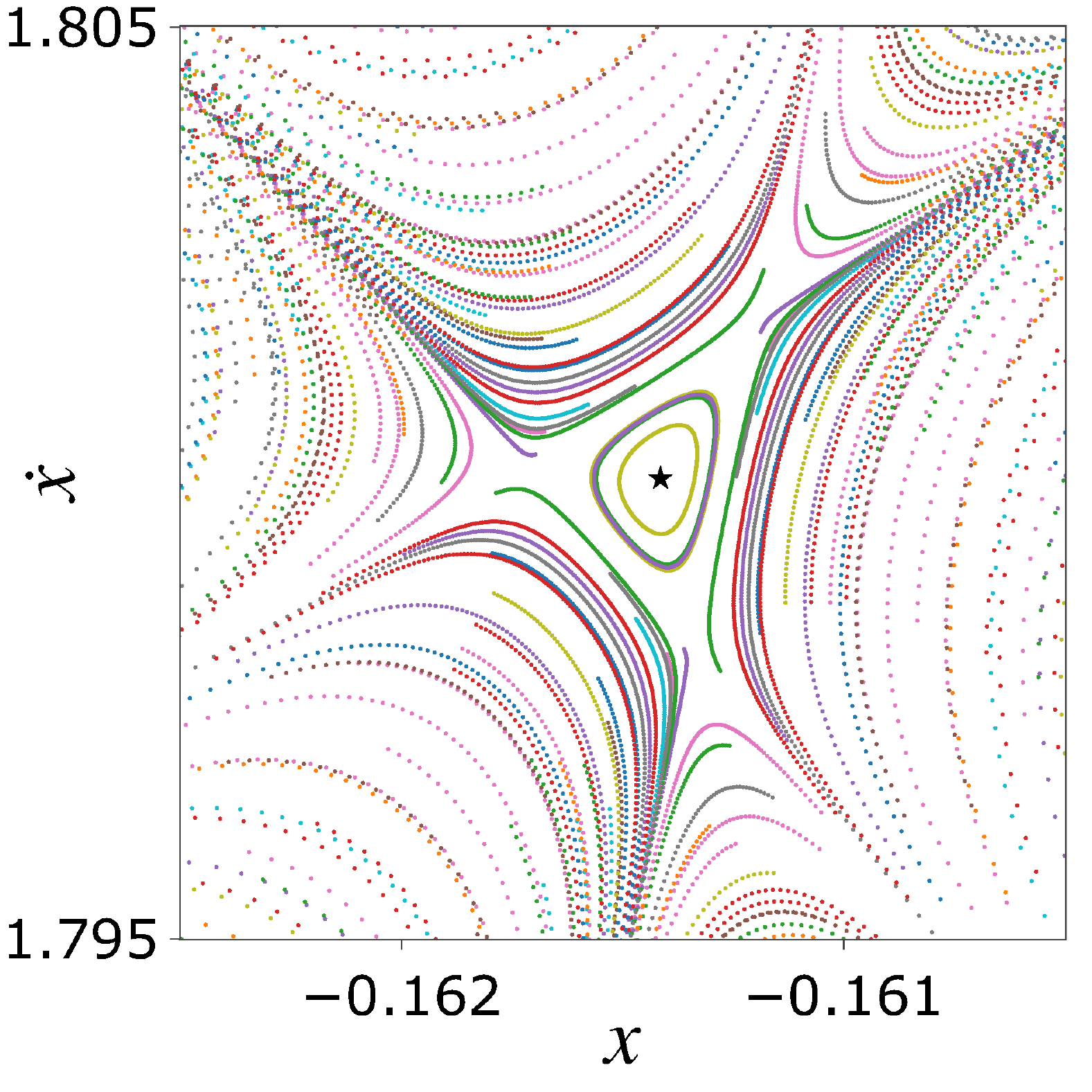}
\caption{$R = 0.895$}
\end{subfigure}
\caption{\small Passing resonance $3 : 1$. When the distance between the walls equals to the critical value $R^* \approx 0.8932$, three separatrices intersect at a fixed point corresponding to a periodic solution without sticking. As $R$ is varied, the separatrices move away from the fixed point forming a triangle which bounds the island of invariant curves surrounding the fixed point. Thus, the fixed point is neutrally stable when $R$ is close to $R^*$ and unstable when $R = R^*$. Parameters are $\omega = 2\pi$, $F = \omega^2$, $f = 24.65$.\label{fig:resonance}}
\end{figure}

In particular, at least a small invariant region $\Omega$ of Hamiltonian dynamics exists around every fixed point  of non-resonance center type because  such a fixed point of $\Phi$ is surrounded by invariant curves, and $\Phi$ preserves the area within the invariant domain bounded by an invariant curve (see Figure~\ref{fig:resonance}). A stable fixed point of center type corresponds to a stable non-sticking periodic solution. For example, symmetric periodic solutions $x(t)=-x(t+\pi/\omega) + r + l$ are defined by the equations
\begin{equation}\label{nonstick}
x_{1, 2}(t) = \begin{cases}
\begin{aligned}
-\frac{F} {\omega^{2}} \cos{(\omega t)} - &\frac f 2 \left(t - \frac{\psi_{1, 2}}{\omega}\right)^2  + C_{1, 2} \left(t - \frac{\psi_{1, 2}}{\omega}\right) + D_{1, 2}, \\ &  0 \le \omega t < \tau_{1,2},
\end{aligned} \\
\begin{aligned}
-\frac{F} {\omega^{2}} \cos{(\omega t)} + &\frac f 2 \left(t - \frac{\psi_{1, 2} + \pi}{\omega}\right)^2 - C_{1, 2} \left(t - \frac{\psi_{1, 2} + \pi}{\omega}\right) - D_{1, 2}, \\ &  \tau_{1,2} \le \omega t < \pi,
\end{aligned}
\end{cases}
\end{equation} where
\begin{equation*}
\psi_{1} = \pi + \arcsin{\left(\frac{\pi f}{2 F}\right)}, \quad \psi_2 = -\arcsin{\left(\frac{\pi f}{2 F}\right)}, \quad \tau_{1, 2} = \psi_{1, 2} \mp \pi,
\end{equation*}
\begin{equation*}
\quad C_{1, 2} = \frac{f\pi^2 \pm 2 \sqrt{4 F^2 - f^2\pi^2}+2R\omega^2}{2 \pi \omega}, \quad D_{1, 2} = l \mp \frac{\sqrt{4 F^2 - f^2 \pi^2}}{2 \omega^2}.
\end{equation*}
One necessary condition for the existence of solutions $x_{1, 2}(t)$ is evident:
\begin{equation}
\frac{f}{F} \le \frac{2}{\pi}.
\end{equation}
Another condition comes from the fact that $x_{1,2}(t)$ are non-sticking periodic solutions. Note that it is sufficient to guarantee that the minimum of the velocity during the free flight from $x=l$ to $x=r$ is positive. Therefore, the non-sticking condition can be written in a closed form as follows:
\begin{equation}\label{nonstick:cond}
\sqrt{4 F^2-\pi ^2 f^2}-\pi  \sqrt{F^2-f^2}+\pi  f \left(\arcsin\left(\frac{\pi
   f}{2 F}\right)-\arcsin\left(\frac{f}{F}\right)\right)+R \omega ^2>0.
\end{equation}
In particular, when
\begin{equation}\label{fold:bif}
\frac f F = \frac 2 \pi
\end{equation} solutions $x_{1,2}(t)$ coincide, which corresponds to the fold bifurcation. At the bifurcation point, condition~\eqref{nonstick:cond} reduces to
\begin{equation}
\frac{ \omega^{2} R}{F} >\left(\sqrt{\pi^2-4}-2 \arccos{\left(\frac 2 \pi\right)}\right).
\end{equation}

Furthermore, we can conclude that the saddle-center bifurcation of periodic solutions is a generic mechanism generating  invariant regions of Hamiltonian dynamics when the kinetic friction parameter $f$ decreases (see Figures~\ref{fig:many:wide}--\ref{fig:many:2pi}).

Stability and type of a fixed point and the corresponding periodic solution can be determined using the linearization of the map $\Phi$ as follows. Any trajectory of~\eqref{main}--\eqref{impact} can be represented as a sequence of motions and events of the following types: {\it free flight} ($\dot x \neq 0$), {\it sticking} ($\dot x = 0$ over a nonzero interval of time), {\it reflection from the wall} ($x = r, l$) and {\it a turning point} (change of sign of $\dot x$). The Jacobi matrix $\Phi^\prime$ can be obtained as a product of matrices of the corresponding four types listed in Table~\ref{tab:matrices}, with the order of matrices corresponding to the order of motions and events during the time interval $0\le t \le T$ (the matrix corresponding to an instantaneous event is known as {\it saltation matrix}).
\begingroup
\begin{table}[H]
\centering
\caption{Matrices corresponding to each type of motion/event.\label{tab:matrices}}
%{\setlength{\tabcolsep}{10pt}
%\renewcommand{\arraystretch}{1.5}
\begin{ruledtabular}
\begin{tabular}{ll}
Type of motion/event & Matrix \\ \hline
Free flight on  $[t_1, t_2]$ & $\begin{pmatrix}
1 & t_2 - t_1\\
0 & 1
\end{pmatrix}$ \\ %\hline
Sticking &$\begin{pmatrix}
1 & 0 \\
0 & 0
\end{pmatrix}$ \\ %\hline
Reflection off the wall at time $t = t^*$ & $\begin{pmatrix}
-1 & 0 \\
\frac{2 F \cos{(\omega t^*)}}{\dot{x}(t^*)} & -1
\end{pmatrix}$ \\ %\hline
Turning point at time $t=\tilde{t}$ & $\begin{pmatrix}
1 & 0 \\
0 & \frac{|F \cos{(\omega \tilde{t})}|-f}{|F \cos{(\omega \tilde{t})}|+f}
\end{pmatrix}$ \\ %\hline
\end{tabular}
\end{ruledtabular}
%}
\end{table}
\endgroup

In particular, for non-sticking trajectories, $\Phi^\prime$ is the product of matrices of the first and the third type (see Table~\ref{tab:matrices}), hence $\det{(\Phi^\prime)} = 1$, which corresponds to the area preservation. On the other hand, a turning point results in shrinking of the phase area due to $\mid\det{(\Phi^\prime)}\mid < 1$. Further, any interval of sticking implies $\det{(\Phi^\prime)} = 0$ and backward non-uniqueness. For example, this is illustrated by Figure~\ref{fig:narrow:0005a} where the fixed point corresponding to a non-sticking periodic trajectory is embedded into the region of Hamiltonian dynamics while the attracting fixed point and period 5 orbit correspond to
exponentially stable periodic trajectories with turning points. Similarly, the attractor on Figure~\ref{fig:many:wide} corresponds to periodic orbits with two intervals of sticking per period.

In Figure~\ref{fig:regs} we explore the global decomposition of the state space into the
domain $\Omega_{pres}$ within which the time $T$ map $\Phi$ is area preserving ($\mid\det(\Phi^\prime)\mid = 1$)
and the complementary domain $\Omega_{dissip}$ where $\Phi$ is contracting ($\mid\det(\Phi^\prime)\mid <1$)\footnote{The derivative $\Phi'$ can be undefined on certain solutions such as trajectories exhibiting grazing or chatter at the wall.}.
Trajectories starting from $\Omega_{pres}$ have non-zero velocity during the whole time interval $0\le t\le T$ while the trajectories starting from $\Omega_{dissip}$ reach zero velocity at least once during the same interval of time.
%In order to achieve thorough understanding of the described phenomenon, it is important to discuss the global dynamics of the map $\Phi$. As we established earlier, the map $\Phi$ preserves the phase area near the fixed points corresponding to the non-stopping trajectories and shrinks it on the trajectories which contain sticking or turning points. Therefore, it is reasonable to consider a global region $\Sigma$ of the phase space of $\Phi$ within which the area is preserved ($\det(\Phi^\prime) = 1$), in contrast to the region $\Delta$ with $\det(\Phi^\prime) < 1$ which we will call {\it dissipative}.
In particular, $\Omega_{pres}$ contains all invariant islands of Hamiltonian dynamics and $\Omega_{dissip}$ contains all the attractors of the map $\Phi$.
%It is necessary to note that here we consider only those points of the phase space for which the linearization $\Phi^\prime$ is properly defined. For instance, if the point corresponds to the orbit with chattering, then $\Phi^\prime$ at this point does not exist. We have encountered chattering while performing numerical simulations of the system~\eqref{main}--\eqref{impact}.
%Figure~\ref{fig:regions} demonstrates the regions $\Sigma$ and $\Delta$ for a particular parameter set.
%We observe that in general $\Sigma$ is not invariant {\it w.r.t.} the map $\Phi$. For example, trajectory denoted by black dots on Figure~\ref{fig:regions} starts in $\Sigma$ near a saddle fixed point and ``walks into" the dissipative region $\Delta$ in a finite number of iterations. However, islands of Hamiltonian dynamics are, of course, invariant. Furthermore, numerical simulations suggest that dissipative region $\Delta$ is invariant, {\it i.e.} $\Phi(\Delta)\subseteq \Delta$, see Figure~\ref{fig:regs}.
In Figure~\ref{fig:regs}, the set  $\Omega_{dissip}$ is invariant for the map $\Phi$, {\em i.e.}
%In other words,
if a solution of system~\eqref{main}--\eqref{impact} has zero velocity at some point during one period, it also has zero velocity at least once per each subsequent period.
%Yet, it is just a conjecture based on the results of numerical experiments, which requires additional verification and rigorous proof. However, it is beyond the scope of this paper, and that's why, we limit ourselves only to the presentation of numerical evidence.
The set $\Omega_{pres}$ splits into an invariant domain of Hamiltonian dynamics  and
a complementary set of points which under iterations of the map $\Phi$ eventually enter the domain $\Omega_{dissip}$ (a trajectory starting from such point has zero velocity at some future moment).
In Figure~\ref{fig:fig} the invariant part of the domain $\Omega_{pres}$ is bounded by the homoclinic tangle of the saddle fixed point corresponding to the non-sticking unstable periodic solution of fundamental frequency $\omega$.

\section{\label{sec:extensions}Possible extensions and concluding remarks}
The analysis presented above can be equally applied to systems with more complex forcing. For example, let us consider the following modification of system~\eqref{main}:
\begin{equation}\label{xpot}
\ddot{x} +f \,\text{sgn}(\dot x)= F \cos{\left(\frac{\pi x} {2}\right)} \cos{(\omega t )}
\end{equation}
with an $x$-dependent force that vanishes on the walls placed at $x = \pm 1$. %It is easy to see that
Unlike~\eqref{main}, system~\eqref{xpot} has a set $E$ of stable equilibrium points, which are located close to the walls.
Namely, $E$ consists of points with zero velocity at the positions $-1\le x\le -1+\eta$ and $1-\eta\le x\le 1$  with
%solutions with $x(t) = x^\dagger$ and $\dot x(t) = 0$ for all $t$, where the values of $x^\dagger$ are located at the vicinity of the walls. Indeed, if
%\begin{equation}\label{steady:cond}
%x^\dagger  \in E := \left[-1, \; -\eta\right]\cup\left[\eta,\; 1\right]
%\end{equation} where
\begin{equation*}
\eta = \frac 2 \pi \arccos{\left(\frac f F\right)}
\end{equation*}
because at these locations friction exceeds the external forcing.
%then the external forcing $\left|F \cos{(\pi x / 2)} \cos{(\omega t)}\right|$ will never exceed friction $f$, hence, the particle will stay at $x = x^\dagger$. Moreover, if $x \in \text{int}(E)$ then, obviously, the steady state solution $\left(x,\; 0\right)$ is stable, in a sense that, any sufficiently small perturbation of it converges in finite time to $\left(y, \; 0\right)$ with $y \in E$.

\begin{figure}[H]
\begin{subfigure}[t]{0.45\textwidth}
\includegraphics*[width=\columnwidth]{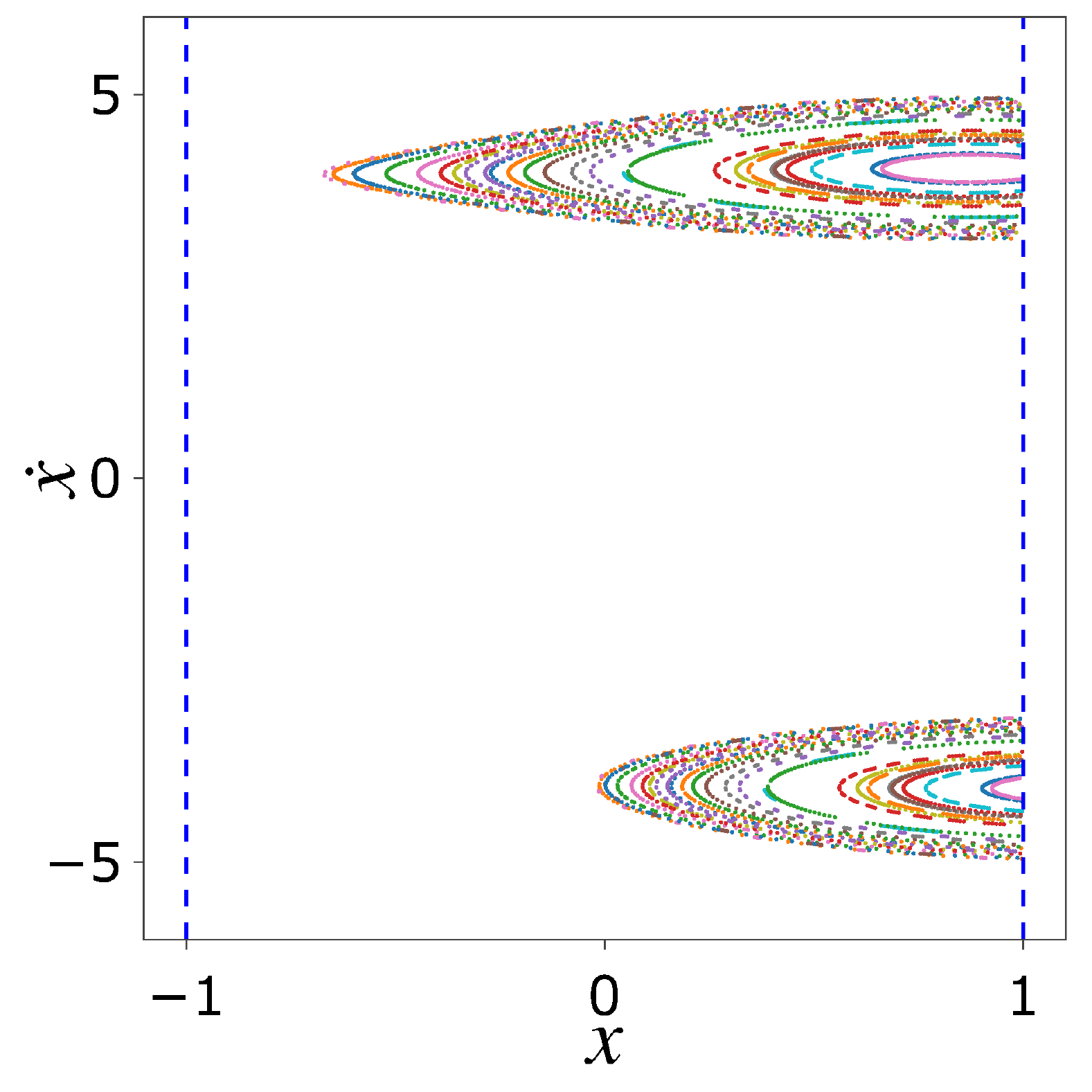}
\caption{\label{fig:xpothamilton}}
\end{subfigure}
\hfill
\begin{subfigure}[t]{0.45\textwidth}
\includegraphics*[width=\columnwidth]{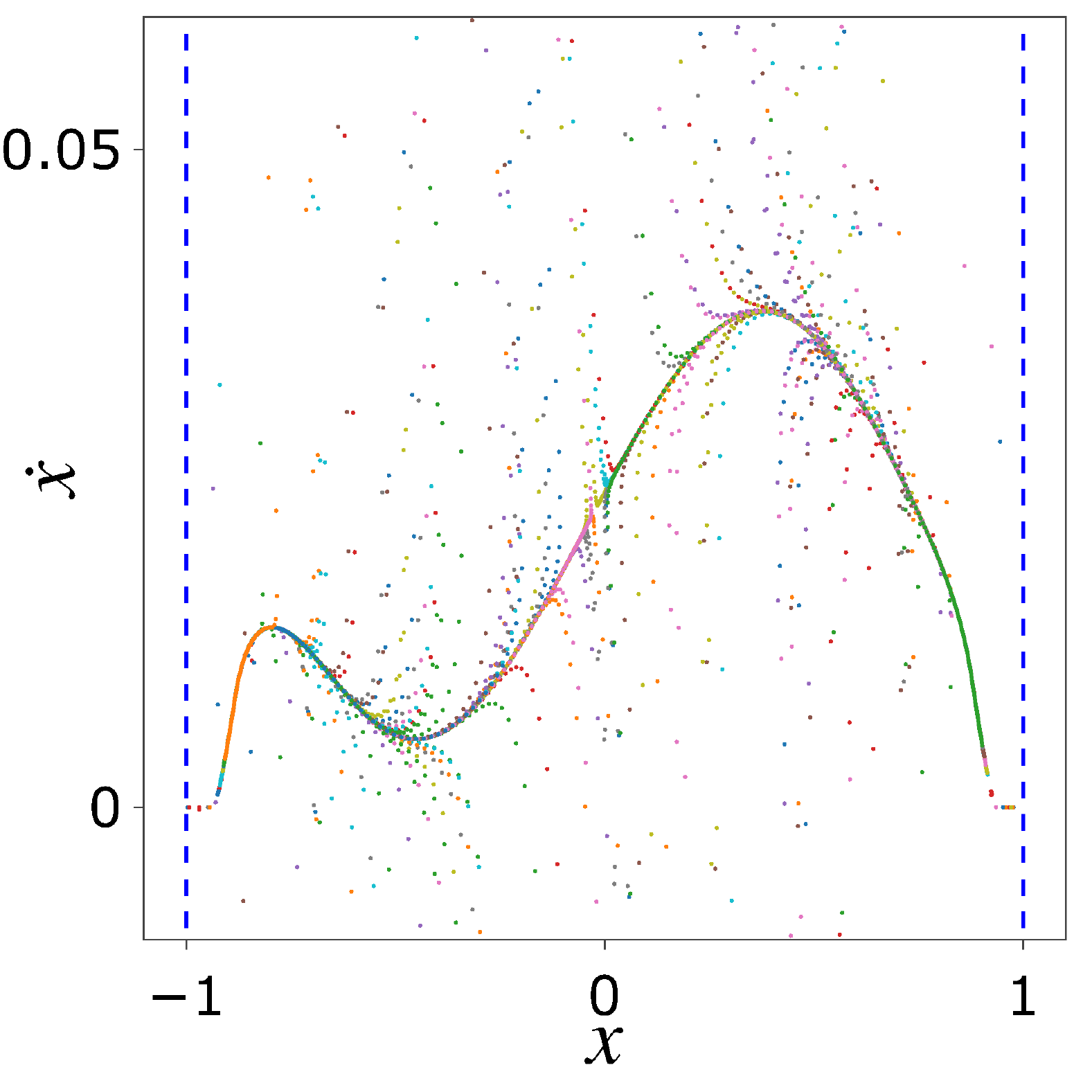}
\caption{\label{fig:xpotnon}}
\end{subfigure}
\caption{\small Orbits of the time $T$ map of system~\eqref{xpot}. Panel $(a)$: Invariant island of Hamiltonian dynamics. Panel $(b)$: Trajectories attracted to steady-state solutions. Parameters are $F = 1$, $f = 0.1$, $\omega = 2\pi$.} \label{fig:xpotno}
\end{figure}

Figure~\ref{fig:xpotno} presents an example of the phase portrait of the time $T$ map for system~\eqref{xpot}. The
set $E$ of equilibrium points is an attractor which co-exists with
 % We see as before it also possesses
 a region of Hamiltonian dynamics with invariant curves surrounding a periodic non-sticking orbit (see Figure~\ref{fig:xpothamilton}). The two components of the attractor $E$ at the walls are connected by an  attractive invariant curve, see Figure~\ref{fig:xpotnon}. Numerical simulations suggest that $E$ is the only attractor.
 In particular, there are no periodic trajectories with stiction or turning points such as in Figures~\ref{fig:narrow-with-friction} and~\ref{fig:many:wide}.

To conclude, %it is important to mention that all main 
a combination of several important
features of the model is crucial for observation of the mixed dynamics. Impact constraints are necessary to switch the velocity direction without passage through zero. Impacts should be elastic, and viscous friction should be absent, to avoid the contraction of the stability map for the non-sticking trajectories. At the same time, the friction coefficient can depend on the particle coordinate.

Further extensions could include multi-particle systems. The simplest case is represented by $N$ identical indistinguishable particles, each satisfying system~\eqref{main}--\eqref{impact}, all placed between the same two walls. If we assume ideal collisions between the particles, then the system is simply equivalent to $N$ independent systems~\eqref{main}--\eqref{impact}. It would be interesting to include interactions between these systems --- for instance, to consider particles of different masses. This is a subject of future work.

% If in two-column mode, this environment will change to single-column format so that long equations can be displayed. 
% Use only when necessary.
%\begin{widetext}
%$$\mbox{put long equation here}$$
%\end{widetext}

% Figures should be put into the text as floats. 
% Use the graphics or graphicx packages (distributed with LaTeX2e).
% See the LaTeX Graphics Companion by Michel Goosens, Sebastian Rahtz, and Frank Mittelbach for examples. 
%
% Here is an example of the general form of a figure:
% Fill in the caption in the braces of the \caption{} command. 
% Put the label that you will use with~\ref{} command in the braces of the \label{} command.
%
% \begin{figure}
% \includegraphics{}%
% \caption{\label{}}%
% \end{figure}

% Tables may be be put in the text as floats.
% Here is an example of the general form of a table:
% Fill in the caption in the braces of the \caption{} command. Put the label
% that you will use with~\ref{} command in the braces of the \label{} command.
% Insert the column specifiers (l, r, c, d, etc.) in the empty braces of the
% \begin{tabular}{} command.
%
% \begin{table}
% \caption{\label{} }
% \begin{tabular}{}
% \end{tabular}
% \end{table}

% If you have acknowledgments, this puts in the proper section head.
\begin{acknowledgments}
The authors are very grateful to Israel Science Foundation (grant 1696/17) for financial support of this work. 
\end{acknowledgments}

% Create the reference section using BibTeX:
\bibliography{biblio}

\end{document}